\documentclass[12pt]{amsart} 
\usepackage{amscd,amsfonts,amsmath,amsxtra,amssymb}
\usepackage[all]{xy}

\newtheorem{sub}{}[section]
\newtheorem{subsub}{}[sub]

\newcommand{\C}{{\mathbb C}}

\renewcommand{\P}{{\mathbb P}}

\newcommand{\F}{{\mathbb F}}

\newcommand{\G}{{\mathbb G}}

\def\hom{\text{hom}}
\def\Hom{\text{Hom}}
\def\mult{\text{mult}}
\def\ext{\text{ext}}
\def\Ext{\text{Ext}}
\def\rk{\text{rk}}
\def\id{\text{id}}
\def\ev{\text{ev}}
\def\Ker{\text{Ker}}
\def\Im{\text{Im}}
\def\GL{\text{GL}}
\def\Stab{\text{Stab}}
\def\Aut{\text{Aut}}
\def\res{\text{res}}

\def\coker{\text{coker}}
\def\length{\text{length}}
\def\span{\text{span}}
\def\codim{\text{codim}}

\newcommand{\kc}{{\mathcal C}}

\newcommand{\ke}{{\mathcal E}}
\newcommand{\kf}{{\mathcal F}}
\newcommand{\kg}{{\mathcal G}}

\newcommand{\ki}{{\mathcal I}}

\newcommand{\kk}{{\mathcal K}}

\newcommand{\ko}{{\mathcal O}}

\begin{document}

\title[Irreducibility and Smoothness]{Irreducibility and Smoothness of the
  moduli space\\ of mathematical 5--instantons over $\mathbf \P_3$}

\author[Coanda]{I.~Coand\u{a}$^{1}$} \address{
  Institute of Mathematics\\
  of the Romanian Academy\\
  P.O.~Box 1-764
\newline RO--70700 Bucharest, Romania}
\email{Iustin.Coanda@imar.ro} 
\author[Tikhomirov]{\;A.~Tikhomirov$^{2}$}
\address{
  Department of Mathematics\\
  State Pedagogical University\\
  Respublikanskaya Str. 108
\newline 150 000 Yaroslavl, Russia}
\email{alexandr@tikhomir.yaroslavl.su} 
\author[Trautmann]{G.~Trautmann$^{3}$}
\address{
  Universit\"at Kaiserslautern\\
  Fachbereich Mathematik\\
  Erwin-Schr\"odinger-Stra{\ss}e
\newline D-67663 Kaiserslautern}
\email{trm@mathematik.uni-kl.de} 
\footnotetext[1]{ Research partially
  supported by the Romanian Academy Grant No.\ 60/2001 and by the\\
\phantom{MA}EURROMMAT
  programme sponsored by the European Commission.}
\footnotetext[2]{Research partially supported by the INTAS programme.}
\footnotetext[3]{Research partially supported by the DFG.}

\begin{abstract}
We prove that the space of mathematical instantons with second Chern class 5
over $\P_3$ is smooth and irreducible. Unified and simple proofs for the same 
statements in case of second Chern class $\leq 4$ are contained.
\end{abstract}

\maketitle
\thispagestyle{empty}

\tableofcontents

\section*{Introduction}

A mathematical instanton bundle is an algebraic vector bundle $\ke$ over
$\P_3(k)$, $k$ an algebraically closed field of characteristic $0$, if it has
rank 2, Chern classes $c_1 = 0$, $c_2 = n > 0$, and if it satisfies the
vanishing conditions $H^0 \ke= 0$ and $H^1 \ke(-2) = 0$.   The name was chosen
in twistor theory in the 1970's, when holomorphic bundles on $\P_3 (\C)$ with
the same properties were considered as counterparts of (anti--)self--dual
Yang--Mills fields on the 4--sphere, see \cite{AHDM}, \cite{Verd},
\cite{Hartsh} for reference.   Let a mathematical instanton bundle with $c_2 =
n$ be called $n$--instanton or $(n,2)$--instanton for short.   The isomorphism
classes of $n$--instantons are the closed points of a coarse moduli scheme
$MI(n)$.   Since its first consideration, it is an open problem whether
$MI(n)$ is smooth and irreducible for any $n$.   An affirmative answer for $n
\le 4$ has been given in several papers, for each $n$ separately, and recently
Katsylo and Ottaviani proved smoothness for $n = 5$, see historical remarks in 
\ref{conj}.

The main result of this paper is that also $MI(5)$ is irreducible, with a new
proof of smoothness included.   The method used also enables a simple and
unified proof for all the previous cases $n \le 4$, see Section
\ref{section5}.

It is well--known that any $n$--instanton is the cohomology of a short complex
$H_n \otimes \ko(-1) \to N \otimes \ko \to H^\ast_n \otimes \ko(1)$, where
$H_n$ and $N$ are $k$--vector spaces of dimensions $n$ and $2n+2$,
respectively, also called the Horrocks construction.   We consider also higher
rank instanton bundles, called $(n,r)$--instantons, which can be constructed
from the same type of complexes with the same $H_n$ but with $\dim N = 2n+r$,
$2 \le r \le 2n$.   We prove that, given an $(n,r)$--instanton $\ke$, one can
choose a linear form $\xi \in H^\ast_n$ such that with $H_{n-1} \cong
\Ker(\xi)$ the induced complex $H_{n-1} \otimes \ko(-1) \to N \otimes \ko \to
H^\ast_{n-1} \otimes \ko(1)$ defines an $(n-1, r+2)$--instanton $\ke_\xi$,
together with a complex $\ko(-1) \to \ke_\xi \to \ko(1)$, whose cohomology is
the original $\ke$, see Section \ref{section3}.   Together with technical
details, this observation enables us to perform induction steps $(n-1, r+2)
\rightsquigarrow (n,r)$ for $n \le 5$ for irreducibility and smoothness.
These induction steps are short for $n \le 4$, see Section \ref{section5}, 
while the
case $(4,4) \rightsquigarrow (5,2)$ is more elaborate.   In Section
\ref{section3}, the induction step is explained in more details.   The proof
of smoothness in each case is achieved by proving that $H^2 S^2 \ke_\xi = 0$
implies $H^2 S^2 \ke = 0$ for suitable $\xi$.    

\section*{Notations}
\begin{itemize}
\item Throughout the paper $k$ will be an algebraically closed field of
  characteristic zero
\item If $E$ is a finite dimensional $k$--vector space, $\P E$ 
  will denote the projective space of $1$--dimensional and $G_mE$ 
  the Grassmannian of $m$--dimensional subspaces.
\item The invertible sheaf of degree $d$ over $\P E$ is denoted
  $\ko_{\P E}(d)$ such that $E^\ast=H^0(\P E, \ko_{\P E}(1))$.
\item $\P=\P V\cong\P_3$ will denote the projective $3$--space over $k$ for a
  fixed vector space $V$ of dimension $4$. We will omit the index $\P$ at the 
  structure sheaf $\ko=\ko_\P$ and at the invertible sheaves
  $\ko(d)=\ko_\P(d)$. 
\item For an $\ko_\P$--module $\kf$ we use the abbreviations
  $\kf(d)=\kf\otimes\ko_\P(d),\ H^i\kf=H^i(\P,\kf)$ and $h^i\kf(d)=\dim_k
  H^i\kf(d)$, and $\ext^i(\kf, \kg)=\dim_k\Ext^i(\kf,\kg)$.
\item If $E$ is a finite dimensional $k$--vector space, the sheaf of sections
  of the trivial vector bundle over a scheme $X$ will be denoted $E\otimes
  \ko_X=E\otimes_k\ko_X$ and $E\otimes\kf$ is written for
  $(E\otimes\ko_X)\otimes_{\ko_X} \kf$.
\item If $\kf$ respectively $\kg$ are coherent sheaves on schemes $X$ resp.\ $Y$,
  $\kf\ \boxtimes\ \kg$ denotes the sheaf $p^\ast\kf\otimes q^\ast\kg$ on
  $X\times_k Y$
  where $p$ respectively $q$ are the projections.
\item If $\kf$ is a coherent sheaf on the scheme $X$ and $Y$ a closed
  subscheme of $X$ we simply write $\kf_Y$ for the restriction on
  $\kf\otimes_{\ko X}\ko_Y=\kf|Y$.
\item Throughout the paper a vector bundle is a locally free sheaf of finite
  rank. If necessary, we will refer specifically to its bundle space as a
  fibration of vector spaces.
\item The Chern classes $c_i(\kf)$ of a coherent sheaf $\kf$ over
  $\P\cong\P_3$ are considered as integers and we also call the triple
  $(c_1(\kf), c_2(\kf), c_3(\kf))$ the Chern class of $\kf$.
\end{itemize}
\vskip1cm

\section{Instanton bundles}\label{section1}

\begin{sub}\rm Let $M(2;0,n,0)$ denote the moduli space of semistable coherent
  sheaves on $\P_3$ of rank $2$ with Chern class $(c_1, c_2, c_3)=(0,n,0)$,
  which had been constructed by M.~Maruyama, \cite{Maru1}, \cite{Maru2}. It
  contains the open set $M^b(2;0,n)$ of stable rank 2 vector bundles on $\P_3$
  with Chern class $(c_1, c_2)=(0,n)$. Recall that a rank $2$ vector bundle
  $\ke$ on $\P_3$ with Chern class $c_1=0$ is stable if and only if
  $H^0\ke=0$. Then it is also simple, i.e.\ $\hom(\ke,\ke)=1$, see
  \cite{OSS}. In this case the Riemann-Roch formula becomes
\[
ext^1(\ke,\ke)-\ext^2(\ke,\ke)=8n-3,
\]
whereas $\Ext^1(\ke,\ke)$ is isomorphic to the tangent space of $M^b(2;0,n)$
at $[\ke]$, the isomorphism class of $\ke$. For large $n$ the spaces
$M^b(2;0,n)$ have many irreducible components and some of them have a much
bigger dimension than $8n-3$. 
\end{sub}
\vskip5mm

\begin{sub} \label{DefInst}\rm {\bf Definition:} A symplectic mathematical 
instanton bundle
  with second Chern class $n\geq 0$ and of rank $r$, or an 
$(n,r)$--instanton for short, is a
  locally free sheaf $\ke$ over $\P=\P V\cong\P_3$ with the following
  properties
  \begin{enumerate}
  \item [(i)] $\ke$ has Chern class $(0,n,0)$ and $2\leq r=rk(\ke)\leq 2n$
  \item [(ii)] $H^0\ke=0$ and $H^1\ke(-2)=0$
  \item [(iii)] $\ke$ admits a symplectic isomorphism
    $\ke\xrightarrow{\varphi}\ \ke^\ast, \varphi^\ast=- \varphi$.
  \end{enumerate}
\end{sub}

Then $\ke$ must have even rank by (iii). In case $rk(\ke)=2$ condition (iii)
can be dropped because then the non--degenerate pairing $\ke \otimes
\ke \to \Lambda^2\ke\cong\ko$ determines a symplectic form, which then is unique
up to a scalar.

If $rk(\ke)=2$, then $\ke$ is stable by condition (ii). Therefore, the open
part 
\[
MI(n)\subset M^b(2;0,n)
\]
defined by the condition $h^1\ke(-2)=0$ is the set of isomorphism classes of
$(n,2)$--instantons, also called $n$--instantons.
\vskip5mm

\begin{sub}\label{conj}{\bf Conjecture:} \rm $MI(n)$ is smooth and irreducible
  of dimension $8n-3$ for $n\geq 1$.

This conjecture emerged in the late 1970's as $n$--instantons were considered
as counterparts of self--dual Yang--Mills fields on the $4$--sphere, see
\cite{AHDM}, \cite{Verd}, in twistor theory.

The conjecture has been proved for $n\leq 4$. For $n=1$ the space $MI(1)$ is
the complement of the Grassmannian $G(2,V)$ in $\P(\Lambda^2 V)$. For $n=2$
the space $MI(2)$ had been described by R.~Hartshorne in \cite{Hartsh} as a
smooth and irreducible fibration. The case $n=3$ was proved by
G.~Ellingsrud--S.A.~Str{\o}mme in \cite{EllStr}. W.~Barth proved in
\cite{Barth2} that $MI(4)$ is irreducible and J.~LePotier in \cite{LePot} that
$MI(4)$ is smooth. Recently P.I.~Katsylo and G.~Ottaviani proved that also
$MI(5)$ is smooth, see \cite{KO}. In this paper we prove that $MI(5)$ is both
smooth and irreducible. The method also yields simple and unified proofs for
the previous results for $n\leq 4$.
\end{sub}
\vskip5mm

\begin{sub}\label{fures} \rm {\bf Further results on $n$--instantons:}

(1) Any $n$--instanton $\ke$ is stable because $H^0\ke=0$ and
$\rk(\ke)=2$. Then the Grauert--M\"ulich theorem states that $\ke$ has trivial
splitting type, i.e.\ for a general line $L$ in $\P_3$ the restricted bundle
$\ke_L$ is isomorphic to $2\ko_L$, see \cite{Barth1}.

(2) A line $L$ is called a jumping line of the $n$--instanton $\ke$ if
$\ke_L\cong \ko_L(-a)\oplus\ko_L(a)$ with $a\neq 0$, and this number is called
the order of the jumping line. It is an easy consequence of the monad
representation of $\ke$, see \ref{beidia} that the set $J(\ke)$ of all jumping
lines of $\ke$ is a hypersurface of degree $n$ in the Grassmannian $\G$ of
lines in $P_3$. Moreover, $n$ is the highest order possible for a jumping
line.

(3) For any $n$--instanton, $h^0\ke(1)\leq 2$, see \cite{BoeTrm}. The
$n$--instantons with $h^0\ke(1)=2$ are called special 'tHooft bundles. They
can be presented as extensions
\[
0\to 2\ko_{\P_3}(-1)\to\ke\to\ko_Q(-n,1) \to 0,
\]
where $Q$ is a smooth quadric in $\P_3$, see \cite{HiNa}. It was shown
there that $MI(n)$ is smooth along the locus of special 'tHooft bundles, which
has dimension $2n+9$.

(4) A plane $P$ in $\P_3$ is called unstable for an  $n$--instanton $\ke$ if
the restricted bundle $\ke_P$ has sections, otherwise stable. If $\ke$ is
special 'tHooft, the unstable planes form a smooth quadric surface in
$\P_3^\ast$, the dual of $Q$ in (3). In \cite{C1} it was proved that an
$n$--instanton $\ke$ is already special 'tHooft if its variety of unstable
planes has dimension $\geq 2$.

(5) In \cite{NueTrm} it was shown that $MI(n)$ is smooth at any $\ke$ with
$h^0\ke(1)\neq 0$ and that $\Ext^2(\ke,\ke(-1))=0$ for such a bundle. The
locus of these bundles has dimension $5n+4$.

(6) P.~Rao, \cite{Rao}, and M.~Skiti, \cite{Skiti}, proved independently that
$MI(n)$ is even smooth along the locus of bundles $\ke$ which admit jumping
lines of highest order $n$. These bundles form a subvariety of dimension
$6n+2$. Moreover, in \cite{Rao} it is proved that $MI(5)$ is smooth at bundles
which have jumping lines of order $4$.

(7) L.~Costa and G.~Ottaviani, \cite{CO}, proved that $MI(5)$ is an affine
scheme by describing the non-degeneracy condition for the monads, see 
\ref{monc}, as the non-vanishing of a hyperdeterminant.

\end{sub}
\vskip5mm
\vskip1cm

\section{Instanton Bundles and Multilinear Algebra}\label{section2}

In this section we describe the monad construction of instantons from
hypernets of quadrics. Throughout this paper
\[
H=H_n 
\]
denotes an $n$-dimensional $k$-vector space,\, $n\geq 1$.
We identify $\Lambda^2(H^\ast\otimes V^\ast)$ with the space of anti--selfdual
$k$--linear maps
\[
H\otimes V\longrightarrow H^\ast\otimes V^\ast.
\]
Each $\omega\in\Lambda^2(H^\ast\otimes V^\ast)$ gives rise to a diagram
\[
\xymatrix{H\otimes V\ar[r]^\omega\ar[d]^{u^\ast} & H^\ast\otimes V^\ast\ar[r] &
  Q\ar[r] & 0\\
N^\ast\ar[r]^\varphi_\sim & N\ar@{>->}[u]_u & & }
\]
where $N=N_\omega$ is the image and $Q=Q_\omega$ the cokernel of $\omega$, and
$\varphi=\varphi_\omega$ is the canonically induced symplectic
isomorphism. Note that the rank of $\omega$ is always even. $\omega$ is called
{\bf non--degenerate} if $\omega(h\otimes v)\neq 0$ for any non--zero
decomposable tensor $h\otimes v$ in $H\otimes V$.
\vskip5mm

\begin{sub}\label{rkstr}\rm{\bf Rank Stratification}: The space
  $\Lambda^2(H^\ast\otimes V^\ast)$ comes with the rank stratification
\[
\Lambda^2(H^\ast\otimes
V^\ast)=\Omega_{4n}\supset\Omega_{4n-2}\supset\ldots\supset\Omega_2, 
\]
where $\Omega_{2m}=\Omega_{2m}(H)=\{\omega\in\Lambda^2(H^\ast\otimes V^\ast)\
|\ \rk(\omega)\leq 2m\}$. This had been considered already in \cite{Tyu} by
A.~Tyurin.

It is easy to prove by standard arguments that each $\Omega_{2m}$ is
irreducible and smooth outside $\Omega_{2m-2}$ of codimension
$\binom{4n-2m}{2}$. The tangent space at $\omega\in
\Omega_{2m}\smallsetminus\Omega_{2m-2}$ is the kernel in
\[
0\to T_\omega\Omega_{2m}\to \Lambda^2(H^\ast\otimes V^\ast)\to \Lambda^2 Q\to
0
\]
which gives the dimension formula. There is the canonical decomposition
\[
\Lambda^2(H^\ast\otimes V^\ast)=(S^2H^\ast\otimes
\Lambda^2V^\ast)\oplus(\Lambda^2H^\ast\otimes S^2 V^\ast),
\]
and on the first summand we have the (induced) rank stratification
\[
S^2H^\ast\otimes \Lambda^2 V^\ast=M_{4n}\supset M_{4n-2}\supset\ldots \supset
M_2
\]
with $M_{2m}=M_{2m}(H)=\Omega_{2m}\cap (S^2H^\ast\otimes \Lambda^2V^\ast)$. We
let
\[
\Delta=\Delta(H)=\{\omega\in S^2H^\ast\otimes \Lambda^2 V^\ast\ |\
\omega\, \text{is degenerate}\}
\]
be the subset of degenerate tensors.

{\bf Remark}: One can prove that $\Delta$ is a closed and irreducible
subvariety of $S^2H^\ast\otimes \Lambda^2V^\ast$ of codimension $2(n-1)$ for
$n\geq 2$ and of codimension $1$ for $n=1$.     
\end{sub}
\vskip5mm

\begin{sub}\label{monc} \rm {\bf The bundles $\ke_\omega$}. To each $\omega\in S^2H^\ast\otimes
  \Lambda^2V^\ast$ one can associate the complex
\[
H\otimes \ko_\P(-1)\xrightarrow{\varphi\circ\alpha^\ast} N\otimes
\ko_\P\xrightarrow{\alpha} H^\ast\otimes \ko_\P(1),
\]
where $\alpha$ is the composition of $N\otimes \ko_\P\xrightarrow{u}
H^\ast\otimes V^\ast\otimes \ko_\P$ and $H^\ast\otimes V^\ast\otimes
\ko_\P\xrightarrow{\id \otimes \ev}H^\ast\otimes \ko_\P(1)$. Then the
following are equivalent:
\begin{enumerate}
\item [(i)]   $\alpha$ is surjective
\item [(ii)]  $\varphi\circ\alpha^\ast$ is a subbundle
\item [(iii)] $\omega$ is non--degenerate.
\end{enumerate}

Clearly (i) and (ii) are equivalent by definition. Now $\alpha$ is surjective
if and only if for any $\langle v\rangle\in \P V$ the induced homomorphism
$N\to H^\ast\otimes \langle v\rangle^\ast$ on the fibre is surjective or
equivalently,  that
$H\otimes \langle v\rangle \xrightarrow{\omega} H^\ast\otimes V^\ast$ is
injective, which is the non--degeneracy.

By this observation each $\omega\in M_{2m}\smallsetminus (M_{2m-2}\cup
\Delta),\, n<m\leq 2n,$ gives rise to an associated locally free sheaf
\[
\ke_\omega=\Ker(\alpha)/\Im(\varphi\circ \alpha^\ast).
\]
This bundle is in fact an instanton bundle of rank $2m-2n$ and with Chern
class $(c_1, c_2, c_3)=(0,n,0)$. This follows directly from the defining
complex, also called the {\bf monad} of $\ke_\omega$. Rank and Chern classes
can be computed from those of the sheaves of the complex as well as
$H^0\ke_\omega=0$ and $H^1\ke_\omega(-2)=0$. Because $\ke^\ast_\omega$ is the
cohomology of the dual complex, it follows that the symplectic isomorphism
$N^\ast\xrightarrow{\varphi} N$ induces a symplectic isomorphism
$\ke_\omega^\ast\xrightarrow{j_\omega} \ke_\omega$.

{\bf Remark}: $M_{2n}\subset\Delta$ because otherwise a non--degenerate
$\omega\in M_{2n}\smallsetminus M_{2n-2}$ would define the bundle
$\ke_\omega=0$ with non--trivial second Chern class. Therefore, $2n+2$ is the
lowest possible rank for a tensor $\omega\in S^2H^\ast\otimes
\Lambda^2V^\ast$, which is non-degenerate, and which then defines an instanton 
of rank $2$.

In the following lemma and proposition it is proved that any $(n,
r)$--instanton arises by the above construction and that the isomorphism
classes of $(n, r)$--instantons are in $1:1$ correspondence with the
equivalence classes of the operators $\omega \in M_{2n+r}\smallsetminus
(M_{2n+r-2}\cup\Delta)$. We use the following {\bf notation} for $2\leq r\leq
2n$.
\[
M(n,r):=M_{2n+r}\smallsetminus (M_{2n+r-2}\cup\Delta)
\]
or 
\[
M(n,r):=\{\omega\in S^2H_n^\ast\otimes \Lambda^2V^\ast\ |\ rk(\omega)=2n+r,
\ \ \omega \text{ non--degenerate}\}.
\]

\vskip5mm

\begin{subsub}\label{ident}{\bf Lemma}: For any $\omega\in M(n,r)$, there are isomorphisms
\[
\begin{array}{lll}
H\cong H^2(\ke_\omega\otimes \Omega^3(1)) & N\cong H^1(\ke_\omega\otimes
\Omega^1)& Q\cong H^1\ke_\omega\\
H^\ast\cong H^1\ke_\omega(-1) & N^\ast\cong H^2(\ke_\omega\otimes\Omega^2) &
\end{array}
\]
which are compatible with Serre duality and the symplectic isomorphism
$\ke^\ast_\omega\cong\ke_\omega$, making the following diagram commute
\begin{tiny}
\[
\xymatrix{H\otimes V\ar[r]\ar[d]^\approx & N^\ast\ar[r]^\varphi_\sim \ar[d]^\approx&
  N\ar[r]\ar[d]^\approx & H^\ast\otimes V^\ast\ar[d]^\approx\ar[r] & Q
  \ar[d]^\approx\ar[r] & 0\\
H^2(\ke_\omega\otimes~\Omega^3(1))\otimes V\ar[r] & H^2(\ke_\omega\otimes
  \Omega^2) & H^1(\ke_\omega\otimes\Omega^1)\ar[l]^\sim_\delta \ar[r] &
  H^1(\ke_\omega(-1))\otimes V^\ast\ar[r] & H^1\ke_\omega\ar[r] & 0
}
\]
\end{tiny}
\end{subsub}
\end{sub}

\begin{proof}The homomorphisms of the bottom row of the diagram are induced by
  the natural homomorphism $\Omega^3(1)\otimes V\to \Omega^2$ and the exact
  sequences $0\to \Omega^2\to \Lambda^2 V^\ast\otimes \ko(-2)\to \Omega^1\to
  0$ and $0\to \Omega^1\to V^\ast\otimes \ko(-1)\to \ko\to 0$ of the Koszul
  complex of $V^\ast\otimes \ko(-1)\to \ko$. Then the isomorphisms of the
  lemma and the commutativity of the diagram follow by tracing all data from
  the defining complex of $\ke_\omega$, using the functoriality of
  Serre-duality.   See also \cite{BH} for monads of vector bundles on
  projective spaces. 
\end{proof}

\vskip5mm

\begin{sub}\label{param}{\bf Proposition}: (a) Let $\ke$ be an
  $(n,r)$--instanton with symplectic isomorphism $\ke^\ast\xrightarrow{j} \ke,
  2\leq r\leq 2n$. Then there exists a non--degenerate $\omega\in
  S^2H^\ast\otimes \Lambda^2 V^\ast$ of rank $2n+r$ such that 
$(\ke,j)\cong(\ke_\omega, j_\omega)$.

(b) $(\ke_\omega, j_\omega)\cong(\ke_{\omega'}, j_{\omega'})$ if and only if
$\omega$ and $\omega'$ are in the same $\GL(H)$--orbit in  $S^2H^\ast\otimes \Lambda^2
V^\ast$. 

(c) The stabilizer $\Stab(\omega)$ of $\omega$ in $\GL(H)$ is isomorphic to the
automorphism group $\Aut(\ke_\omega, j_\omega)$.

\begin{proof}(a) By \ref{ident} it is enough to show that the Beilinson
  spectral sequence of $\ke$ results in the complex
\[
0\to H^2(\ke\otimes \Omega^3(1))\otimes \ko(-1)\xrightarrow{\beta}
H^1(\ke\otimes \Omega^1)\otimes \ko\xrightarrow{\alpha} H^1(\ke(-1))\otimes
\ko(1)\to 0,
\]
where $\alpha$ resp.\ $\beta$ correspond to the natural homomorphisms
\[
H^2(\ke\otimes \Omega^3(1))\otimes V\to H^2(\ke\otimes
\Omega^2)\underset{\sim}{\leftrightarrow} H^1(\ke\otimes \Omega^1)
\]
resp.
\[
H^1(\ke\otimes \Omega^1)\to H^1\ke(-1)\otimes V^\ast,
\]
which are Serre--dual to each other. Now for any coherent sheaf $\kf$ on
$\P_n$ there is the Beilinson I complex
\[
0\to C^{-n}(\kf)\to\cdots\to C^0(\kf)\to\cdots\to C^n(\kf)\to 0
\]
with terms $C^p(\kf)=\underset{p=i-j}{\oplus} H^i(\kf\otimes
\Omega^j(j))\otimes \ko(-j)$, which is exact except at $C^0(\kf)$ and which
has $\kf$ as its cohomology at $C^0(\kf)$. In our case for $\kf=\ke(-1)$ this
complex reduces in fact to
\[
0\to H^1(\ke\otimes \Omega^2(1))\otimes \ko(-2)\to
H^1(\ke\otimes\Omega^1)\otimes\ko(-1)\to H^1(\ke(-1))\otimes \ko\to 0
\]
by verifying that the instanton conditions imply the vanishing of the other
terms of the Beilinson complex. Moreover, the Koszul sequence
$0\to\Omega^3\to\Lambda^3 V^\ast\otimes\ko(-3)\to\Omega^2\to 0$ induces the
isomorphism $H^1(\ke\otimes \Omega^2(1))\cong H^2(\ke\otimes\Omega^3(1))$
because $H^i\ke(-2)=0$ for $i=1,2$. Now $\ke$ is the cohomology of the complex
\[
0\to H^2(\ke\otimes\Omega^3(1))\otimes \ko(-1)\xrightarrow{\beta}
H^1(\ke\otimes \Omega^1)\otimes \ko\xrightarrow{\alpha} H^1(\ke(-1))\otimes
\ko(1)\to 0.
\]
It follows from the Riemann--Roch formula that 
\[
\chi\ke(d)=r\binom{d+3}{3}-n(d+2)
\]
and from the instanton conditions that
\[
h^1\ke(-1)=n\quad\text{ and }\quad h^1(\ke\otimes \Omega^1)=2n+r.
\]
Thus $\ke$ defines a non--degenerate $\omega\in S^2H^\ast\otimes \Lambda^2
V^\ast
$
via
\[
H^2(\ke\otimes \Omega^3(1))\otimes V\to H^2(\ke\otimes
\Omega^2)\xleftarrow[\sim]{} H^1(\ke\otimes \Omega^1)\to
H^1\ke(-1)\otimes V^\ast
\]
together with an isomorphism $H^\ast\cong H^1\ke(-1)$, such that $(\ke,
j)\cong(\ke_\omega, j_\omega)$.

(b) and (c) follow from (a) and \ref{ident} by the isomorphisms of the complexes.
\end{proof}
\end{sub}
\vskip5mm

\begin{sub}{\bf Corollary}: If $\ke_\omega$ is simple, then
  $\Stab_{\GL(H)} (\omega)=\{\pm\id_H\}$.
\end{sub}
\vskip5mm

\begin{sub}\label{beidia}{\bf Remark}: \rm There is also a Beilinson II monad for
  an instanton bundle $\ke_\omega$. This can be treated in the same way as the
  above Beilinson I monad. Any $\omega\in S^2H^\ast\otimes \Lambda^2 V^\ast$
  defines an operator $H\to H^\ast\otimes \Lambda^2 V^\ast$ which is symmetric
  with respect to $H$ and the exact sequence $0\to N\xrightarrow{u}
  H^\ast\otimes V^\ast\to Q\to 0$ as above. Now combined with the Koszul homomorphisms
\[
\Lambda^2 V^\ast\otimes \ko(-1)\to \Omega^1(1)\quad \text{ and }\quad
\Omega^1(1)\to V^\ast\otimes \ko
\]
we obtain the complex
\[
0\to H\otimes \ko(-1)\xrightarrow{\mu} H^\ast\otimes
\Omega^1(1)\xrightarrow{\nu} Q\otimes \ko\to 0.
\]
If $\omega$ is non--degenerate, this complex defines the $(n,r)$--instanton
$\ke_\omega$, where $\rk(\omega)=2n+r$. This could be shown as in \ref{ident},
\ref{param}, but follows already from the following commutative diagram induced
by $\omega$, which is a direct transformation between the two monad types.
\[
\xymatrix{ & 0\ar[d] & 0\ar[d] & &\\
0\ar[r] & H\otimes \ko(-1)\ar[d]\ar[r]^\mu & H^\ast\otimes
\Omega^1(1)\ar[d]\ar[r] & Q \otimes \ko\ar[r]\ar@{=}[d] & 0\\
0\ar[r] & N\otimes \ko\ar[d]\ar[r]^-u & H^\ast\otimes V^\ast\otimes
\ko\ar[d]\ar[r] & Q\otimes \ko\ar[r] & 0\\
 & H^\ast\otimes \ko(1)\ar[d] \ar@{=}[r] & H^\ast\otimes \ko(1)\ar[d] & & \\
 & 0 & 0 & &}
\]
\end{sub}
\vskip5mm

\begin{sub}{\bf The morphism $\mathbf
    M_{2n+2}\smallsetminus\Delta\overset{b}{\twoheadrightarrow} MI(n)$} 

\rm According to Proposition \ref{param} we are given a surjective map
$\omega\overset{b}{\mapsto}[\ke_\omega]$ from $M_{2n+2}\smallsetminus\Delta$
to $MI(n)$ whose fibres are the orbits under the action of $\GL(H)$. The map
$b$ is the underlying map of a morphism, because there is a universal monad
over $(M_{2n+2}\smallsetminus\Delta)\times\P$ with the universal family of
$n$--instantons as cohomology. Then $b$ is the modular morphism onto the open
part $MI(n)$ of the Maruyama scheme. It can be shown that $b$ is a geometric
quotient and a principal $GL(H)/\{\pm 1\}$--bundle in the etale topology.
However the latter fact will not be used in this paper.
\end{sub}
\vskip5mm

\begin{sub}\label{expd} \rm {\bf Criterion for irreducibility}

Since $\Omega_{2n+2}$ is irreducible of codimension $\binom{2n-2}{2}$ in
$\Lambda^2(H^\ast\otimes V^\ast)$, it follows that every irreducible component
of $M_{2n+2}$ has codimension $\leq \binom{2n-2}{2}$ in $S^2H^\ast\otimes
\Lambda^2 V^\ast$ or dimension $\geq 3n(n+1)-\binom{2n-2}{2}=n^2+8n-3$. From
this observation we obtain
\vskip5mm

\begin{subsub}\label{crirr} {\bf Lemma}: 
Any component of $MI(n)$ has dimension $\geq 8n-3$ and
$MI(n)$ is irreducible of dimension $8n-3$ if and
  only if $M_{2n+2}\smallsetminus \Delta$ is irreducible of the expected
  codimension $\binom{2n-2}{2}$ in $S^2H^\ast\otimes \Lambda^2 V^\ast$.
\end{subsub}

\rm For the proof use the following lemma.

\vskip5mm

\begin{subsub}\label{cricomp} {\bf Lemma}: Let $X\xrightarrow{f} Y$ be a
  morphism of reduced schemes, let $Y$ be irreducible and assume that each
  fibre $f^{-1}(f(x)),\, x\in X,$ is irreducible of constant dimension $e$, and
  that each irreducible component of $X$ is of dimension $\geq\dim Y+e$. Then
  $X$ is irreducible.

\begin{proof} Let $X_1, \ldots, X_p$ be the irreducible components of
  $X$. According to the theorem of the dimension of the fibres, $f|X_i$ is
  dominant over $Y$ for any $i$. Let $X'_i\subset X_i$ be the complement of
  the other components and let $y\in\cap f(X'_i)$. Then
\[
f^{-1}(y)=\cup (f^{-1}(y)\cap X_i)
\]
and the intersections $f^{-1}(y)\cap X_i$ are distinct closed subsets of
$f^{-1}(y)$. As $f^{-1}(y)$ is irreducible, it follows that $p=1$.
\end{proof}
\end{subsub}
\vskip5mm

\begin{subsub}\label{rmtyu} \rm {\bf Remark}: In \cite{Tyu2}, Example 1,
  A.~Tyurin showed that $M_{2n}$ and (consequently) $M_{2n+2}$ have not the
  expected codimension $\binom{2n}{2}$ respectively $\binom{2n-2}{2}$ in
  $S^2H^\ast\otimes \Lambda^2 V^\ast$.
\end{subsub}
\end{sub}
\vskip5mm

\begin{sub}\rm {\bf Criterion for smoothness}

Concerning the smoothness of $MI(n)$, let $\ke$ be an $n$--instanton (of rank 2). Because
$[\ke]$ is a stable point of $MI(n)$, this is a smooth point if
$\Ext^2(\ke,\ke)=0$. But since $\ext^1(\ke,\ke)-\ext^2(\ke,\ke)=8n-3$ the
vanishing of $\Ext^2(\ke,\ke)$ is also necessary if $[\ke]$ is a smooth point
and $MI(n)$ has dimension $8n-3$ at $[\ke]$. On the other hand, we have
\[
\Ext^2(\ke,\ke)\cong H^2(\ke^\vee\otimes \ke)\cong H^2(\ke\otimes\ke)\cong
H^2(S^2\ke).
\]
Now for any $(n,r)$--instanton we have the following
\vskip5mm

\begin{subsub}\label{crism} {\bf Lemma}: Let $\omega\in S^2H^\ast\otimes
  \Lambda^2 V^\ast$ be non--degenerate of rank $2m,\ n<m\leq 2n$, and let
  $\ke=\ke_\omega$. Then the following conditions are equivalent.
  \begin{enumerate}
  \item[(i)] $H^2(S^2\ke)=0$
  \item[(ii)] $\Omega_{2m}$ and $S^2H^\ast\otimes\Lambda^2 V^\ast$ intersect
    transversally in $\omega$ inside $\Lambda^2(H^\ast\otimes V^\ast)$.
\item [(iii)] $M_{2m}$ is smooth at $\omega$ of expected dimension $5n-5n^2+8mn-2m^2-m$.
  \end{enumerate}

\begin{proof} Transversality in (ii) means that the tangent space
  $T_\omega\Omega_{2m}$ and $S^2H^\ast\otimes \Lambda^2 V^\ast$ $\span\ 
  \Lambda^2(H^\ast\otimes V^\ast)$. Then (ii) and (iii) are equivalent by
  standard dimension counts, because $\Omega_{2m}$ has codimension
  $\binom{4n-2m}{2}$ and is smooth at $\omega$.

 Next, let us recall a general fact. To any bounded complex
  $K^\bullet$ of vector bundles on a variety one can associate its second
  symmetric power $S^2K^\bullet$ by decomposing $K^\bullet\otimes K^\bullet$
  according to the eigenspaces of the canonical involution on $K^\bullet\otimes
  K^\bullet$. If $K^\bullet$ has only three non--zero terms $K^{-1}\to K^0\to K^1$,
  then $S^2K^\bullet$ is isomorphic to 
\[
0\to \Lambda^2 K^{-1}\to K^{-1}\otimes K^0\to S^2 K^0\oplus(K^{-1}\otimes
K^1)\to K^0\otimes K^1\to \Lambda^2K^1\to 0
\]
whose differentials are naturally deduced from those of $K^\bullet$. Moreover, if
$H^i(K^\bullet)=0$ for $i\neq 0$ and $H^0(K^\bullet)=F$, then $H^i(S^2K^\bullet)=0$
for $i\neq 0$ and $H^0(S^2K^\bullet)\cong S^2 F$. Hence, in our case $S^2\ke$ is
the degree $0$ cohomology of the derived complex
\[
\Lambda^2 H\otimes \ko(-2)\to H\otimes N\otimes \ko(-1)\to (S^2 N\oplus
H\otimes H^\ast)\otimes\ko\to N\otimes H^\ast\otimes \ko(1)\to \Lambda^2
H^\ast\otimes \ko(2).
\]
The terms of this monad are cohomologically acyclic. Hence one can compute the
coho\-mology of $S^2\ke$ just by passing to global sections and then taking
cohomology. Particularly, one gets an exact sequence
\[
N\otimes H^\ast\otimes V^\ast\to \Lambda^2 H^\ast\otimes S^2 V^\ast\to
H^2(S^2\ke)\to 0.
\]
Now condition (ii) is equivalent to 
\[
T_\omega\Omega_{2m}+(S^2H^\ast\otimes \Lambda^2 V^\ast)= \Lambda^2(H^
\ast\otimes V^\ast),
\]
noting that $\omega$ is a smooth point of $\Omega_{2m}$. Taking into account
the exact sequences
\[
N\otimes (H^\ast\otimes V^\ast)\to \Lambda^2(H^\ast\otimes V^\ast)\to
\Lambda^2 Q\to 0
\]
and
\[
0\to S^2 H^\ast\otimes \Lambda^2 V^\ast\to \Lambda^2(H^\ast\otimes V^\ast)\to
\Lambda^2 H^\ast\otimes S^2 V^\ast\to 0,
\]
one finds that (ii) is equivalent to the surjectivity of the composed map
\[
N\otimes H^\ast\otimes V^\ast\to \Lambda^2(H^\ast\otimes V^\ast)\to \Lambda^2
H^\ast\otimes S^2 V^\ast,
\]
which is the map whose cokernel is $H^2(S^2\ke)$. This proves the equivalence
of (i) and (ii).
\end{proof}
\end{subsub}
\vskip5mm

Now, by the criteria for irreducibility and smoothness, it is clear that the
conjecture \ref{conj} on $MI(n)$ is equivalent to
\end{sub}
\vskip5mm

\begin{sub}\label{conj2}\rm {\bf Transcribed Conjecture}: $\Omega_{2n+2}$ and
  $S^2H^\ast\otimes \Lambda^2 V^\ast$ intersect transversally inside
  $\Lambda^2(H^\ast\otimes V^\ast)$ along $M_{2n+2}\smallsetminus\Delta$ and
  this intersection is irreducible.

Note, that in this conjecture, with $\Delta$ also $M_{2n}$ has been extracted
from $M_{2n+2}$ because $M_{2n}\subset\Delta$. According to Tyurin's example
the whole of $M_{2n+2}$ has components of excessive dimension for large $n$.
\end{sub}
\vskip5mm

\begin{sub}\rm {\bf Remark}\label{remsm}: Katsylo--Ottaviani gave in \cite{KO}
  the following interpretation of the transversality condition. Dualizing the
  sequence with $H^2(S^2\ke)$ as cokernel in the proof of Lemma \ref{crism},
  we have the exact sequence
\[
0\to H^2(S^2\ke)^\ast\to \Lambda^2 H\otimes S^2 V\to N^\ast\otimes H\otimes V.
\]
Interpreting the elements of $\Lambda^2 H \otimes S^2 V$ as anti--selfdual
linear maps $H^\ast\otimes V^\ast\to H\otimes V$ and the elements of
$N^\ast\otimes H\otimes V$ as linear maps $N\to H\otimes V$, the morphism in
the above sequence can be described by $\sigma\mapsto \sigma\circ u$, where
$u$ is the inclusion of $\Im(\omega)=N\subset H^\ast\otimes V^\ast$. Therefore,
\[
H^2(S^2\ke)^\ast\cong\{\sigma\in \Lambda^2 H\otimes S^2 V\ |\
\sigma\circ\omega=0\}
\]
and the three conditions of \ref{crism} are equivalent to
\[
\text{(iv) if}\ \sigma\circ\omega=0\quad \text{ for } \sigma\in\Lambda^2
H\otimes S^2 V,\ \text{ then } \sigma=0.
\]
In order to illustrate this point of view and for later use, we prove
\end{sub}
\vskip5mm

\begin{sub}\label{hcase} {\bf Proposition}: $M(n, 2n-2) = M_{4n-2}\smallsetminus(\Delta\cup
  M_{4n-4})$ is smooth and (obviously) of codimension $1$ in $S^2H^\ast\otimes
  \Lambda^2V^\ast$ for $n\geq 2$.

\begin{proof} Let $\omega\in M_{4n-2}\smallsetminus(\Delta\cup
  M_{4n-4})$ and $\sigma\in \Lambda^2 H\otimes S^2 V$ such that
  $\sigma\circ\omega=0$. Since both $\omega$ and $\sigma$ are anti--selfdual,
  it follows that $\omega\circ \sigma=0$. Now $\rk(\omega)=4n-2$ implies
  $\rk(\sigma)\leq 2$. One can now easily show (see e.g.\ \cite{Tyu},
  Proposition 2.1.1) that in this case $\sigma$ is decomposable,
  $\sigma=\eta\otimes f$, with $\eta\in\Lambda^2 H$ of rank $\leq 2$ and $f\in
  S^2 V$ of rank $\leq 1$. Consequently, if $\sigma\neq 0$, then
  $\Im(\sigma)\subset H\otimes V$ contains decomposable vectors, contradicting $\omega\not\in\Delta$.
\end{proof}\end{sub}
\vskip1cm

\section{The Method}\label{section3}

Let $\ke$ be an $(n,r)$--instanton and $H^\ast=H^\ast_n\cong H^1\ke_\omega(-1)$ as above
with $\omega\in S^2H^\ast\otimes \Lambda^2 V^\ast\setminus \Delta$ of rank
$2n+r$. To any $\xi\in H^\ast$ with kernel $\bar{H}$ we can associate the
restriction $\bar{\omega}=\omega_\xi=\res_\xi(\omega)$ of $\omega$ to $S^2\bar{H}^\ast\otimes
\Lambda^2 V^\ast$ such that we have the diagram
\[
\xymatrix{
\bar{H}\otimes V\ar[r]^{\bar{\omega}}\ar@{>->}[d] &
  \bar{H}^\ast\otimes V^\ast\\
H\otimes V\ar[r]^\omega & H^\ast\otimes V^\ast\ar@{->>}[u]
}.
\]
Then $\rk(\bar{\omega})\leq\rk(\omega)$.

The choice of $\xi$ gives rise to a diagram
\[
\xymatrix{ & 0\ar[d] & & 0\ar[d]&\\
 & \bar{H}\otimes \ko(-1)\ar[d]\ar[dr]^{\bar{\beta}} & &\ko(1)\ar[d] & \\
0\ar[r] & H\otimes \ko(-1)\ar[d]\ar[r]^-\beta & N\otimes
 \ko\ar[r]^-\alpha\ar[dr]^{\bar{\alpha}} & H^\ast\otimes \ko(1)\ar[r]\ar[d] & 0\\
 & \ko(-1)\ar[d] & & \bar{H}^\ast\otimes\ko(1) \ar[d]& \\
 & 0 & & 0 &}
\]
such that $\bar{\beta}, \bar{\alpha}$ constitute a monad for a bundle
$\bar{\ke}$ with an induced symplectic isomorphism, and such that $\bar{\beta},
\bar{\alpha}$ are induced by $\bar{\omega}$. Then $\bar{\ke}$ has
rank $r+2$ and second Chern class $n-1$. Any splitting homomorphism of
$H\otimes \ko(-1)\to \ko(-1)$ then induces a unique homomorphism
$\ko(-1)\to\bar{\ke}$ and by duality a monad
\[
0\to \ko(-1)\to\bar{\ke}\to\ko(1)\to 0
\]
whose cohomology is again $\ke$. This observation will be used to perform an
induction $(n-1, r+2)\rightsquigarrow (n,r)$ for irreducibility and smoothness
of the spaces of $(n,r)$--instantons for $n\leq 5$.

Immediate relations between the tensors $\omega$ and $\bar{\omega}$ and the
corresponding sheaves $\ke$ and $\bar{\ke}$ are stated in the next two lemmata.
\vskip5mm

\begin{sub}\label{rkl} {\bf Lemma}: Let $\bar{\omega}$ be the restriction of
  $\omega$ for $\xi\in H^\ast$. The following conditions are equivalent:
  \begin{enumerate}
  \item [(i)]   $\rk(\bar{\omega})=\rk(\omega)$
  \item [(ii)]  $\Im(\omega)\cap(\xi\otimes V^\ast)=0$
  \item [(iii)] The multiplication map $H^1\ke(-1)\otimes V^\ast\to H^1\ke$
    restricts to an injection\\$\xi\otimes V^\ast\hookrightarrow H^1\ke$.
  \item [(iv)] $H^0\bar{\ke}=0$
  \end{enumerate}

\begin{proof}We have the diagram
\[
\xymatrix{\bar{H}\otimes V\ar[r]^{\bar{\omega}}\ar[d]^j & \bar{H}^\ast\otimes
  V^\ast & \\
H\otimes V\ar[r]^\omega\ar[d]^{u^\ast} & H^\ast\otimes
  V^\ast\ar[u]_{j^\ast}\ar[r]^-\mu & Q\ar[r] & 0\\
N^\ast\ar[r]^\varphi_\approx & N \ar[u]_u& &
}
\]
where $N$ is the image of $\omega, \Ker(j^\ast)=\xi\otimes V^\ast$, and $\mu$
is isomorphic to the multiplication map $H^1\ke(-1)\otimes V^\ast\to
H^1\ke$. It follows that $j^\ast\circ u$ is injective if and only if the image
of $\bar{\omega}$ has the same dimension as $N$, which proves the equivalence
of (i) and (ii). The other statements are immediately seen to be equivalent to
the injectivity of $j^\ast\circ u$. For (iv), note, that
$N\to\bar{H}^\ast\otimes V^\ast$ corresponds to the right part of the monad of
$\bar{\ke}$.

\end{proof}\end{sub}

If the conditions of the lemma are satisfied, then $\bar{\ke}$ is
  an $(n-1, r+2)$--instanton, provided $r<2n$. In this case we have
\vskip5mm

\begin{sub}\label{vanind} {\bf Lemma}: Let $\omega\in S^2H^\ast\otimes
  \Lambda^2 V^\ast$ be non--degenerate and $\rk(\bar{\omega})=\rk(\omega)$. If
  $H^2(S^2\bar{\ke})=0$ and $H^1\bar{\ke}(1)=0$, then $H^2(S^2\ke)=0$.

\begin{proof} Because $\ke$ is the cohomology of the monad $\ko(-1)
  \to\bar{\ke}\to \ko(1),\, S^2\ke$ is the cohomology of the induced monad, see
  proof of Lemma \ref{crism},
\[
\bar{\ke}(-1)\to\ko\oplus S^2\bar{\ke}\to\bar{\ke}(1).
\]
The computation of $H^2(S^2\ke)$ from this monad gives the vanishing.
\end{proof}
\end{sub}
\vskip5mm

\begin{sub}\label{indproc}\rm {\bf Construction lemma}

For $n\geq 3$ and $2\leq r\leq 2n-2$ and for any $\xi\in H^\ast$ we consider
the open subset
\[
M(n,r)_\xi:=\{\omega\in M(n,r)\ |\ \Im(\omega)\cap(\xi\otimes V^\ast)=0\}.
\]
This is exactly the open subset of $M(n,r)$ which is mapped to $M(n-1, r+2)$
under the restriction map
\[
\res_\xi:\omega\mapsto \omega_\xi=\omega \ |\ S^2\bar{H}^\ast\otimes \Lambda^2
V^\ast.
\]
We also use the notation
\[
M^0(n, r)=\{\omega\in M(n, r)\ |\
H^1\ke_{\omega}(1)=0\}.
\]
Note that this set might be empty if $h^0\ke_{\omega}(1)-h^1\ke_{\omega}(1)
=4r-3n< 0,$
while $h^2\ke_{\omega}(1)=h^3\ke_{\omega}(1)=0$. 
Next we let 
\[
M(n,r)'_\xi=\{\omega\in M(n,r)_\xi\ |\ \res_\xi(\omega)\in M^0(n-1, r+2)\}
\]
be the inverse image. For the fibres of $\res_\xi$ we have the following 
lemma, which enables the induction. Note, that in this lemma the fibre 
$R(\bar{\omega})$ may contain degenerate $\omega$'s.
\vskip5mm

\begin{subsub}\label{fibresxi} {\bf Lemma}: a) For any $\bar{\omega}\in M(n-1,
  r+2)$ there is an isomorphism 
\[
R(\bar{\omega}): = \{\omega\in S^2H_n^\ast\otimes \Lambda^2V^\ast\ |\ 
\res_\xi(\omega) = \bar{\omega},\quad 
\rk(\omega)=\rk(\bar{\omega})\}\cong\Hom_{\ko_\P}(\kc_{\bar{\omega}},
\ko_\P(1))
\]
where $\kc_{\bar{\omega}}$ denotes the cokernel of the left part of the monad
of $\ke_{\bar{\omega}}$, and
\[
\dim\Hom_{\ko_P}(\kc_{\bar{\omega}}, \ko_\P(1))=n-1+h^0\ke_{\bar{\omega}}(1).
\]

b) $\omega\in R(\bar{\omega})$  is non--degenerate if and only if the 
corresponding homomorphism
induces an epimorphism $\ke_{\bar{\omega}}\to \ko_\P(1)$. In that case
$\omega$ defines an $(n,r)$--instanton bundle $\ke_\omega$ which is also the cohomology
of the self--dual monad $\ko_\P(-1)\to\ke_{\bar{\omega}}\to\ko_\P(1)$ defined
by the epimorphism.
\end{subsub}

\begin{proof} a) For $\xi$ fixed we may choose a decomposition $H=H_1\oplus
  \bar{H}$ such that $\xi$ induces an isomorphism $H_1\cong k$. Let
  $\bar{\omega}$ decompose into 
\[
\bar{H}\otimes
V\overset{\bar{u}^\ast}{\twoheadrightarrow}N^\ast\xrightarrow[\approx]{\varphi}N
\overset{\bar{u}}{\rightarrowtail}\bar{H}^\ast\otimes
V^\ast
\]
with associated monad
\[
\bar{H}\otimes \ko_\P(-1)\xrightarrow{\bar{\beta}} N\otimes
\ko_\P\xrightarrow{\bar{\alpha}} \bar{H}^\ast\otimes \ko_\P(1).
\]
With $\kc_{\bar{\omega}}=\coker(\bar{\beta})$ we have the exact sequence
\[
0\to \ke_{\bar{\omega}}\to\kc_{\bar{\omega}}\to \bar{H}^\ast\otimes \ko_\P(1)\to 0.
\]
Any $u_1\in\Hom (N, H_1^\ast\otimes V^\ast)$ gives rise to a skew--symmetric
operator
\[
\omega=\left(
  \begin{array}{l|c}
u_1\circ\varphi\circ u_1^\ast & u_1\circ\varphi\circ \bar{u}^\ast\\
\hline
\bar{u}\circ\varphi\circ u^\ast_1 & \bar{\omega}
  \end{array}\right)
\]
\[
(H_1\oplus\bar{H})\otimes
V\xrightarrow{\omega}(H_1^\ast\oplus\bar{H}^\ast)\otimes V^\ast
\]
which factors also through $\varphi$ by its definition, so that
$\rk(\omega)=\rk(\bar{\omega})$ and $\omega$ is symplectic. However, the
component $u_1 \circ \varphi \circ \bar{u}^\ast : \bar{H} \otimes V \to H_1^\ast
\otimes V^\ast$ is not necessarily skew with respect to $V$.   This is the
case if and only if $\omega \in R(\bar{\omega})$ or if and only if the
composition
\[
\bar{H} \otimes \ko_\P (-1) \xrightarrow{\beta} N \otimes \ko_\P
\xrightarrow{\tilde{u}_1} H_1^\ast \otimes \ko_\P (1)
\]
of the associated sheaf homomorphismus is zero, or, if and only if
$\tilde{u}_1$ factors through $\kc_{\bar{\omega}}$.   Let $\Hom(N, H_1^\ast
\otimes V^\ast)^\prime$ denote the subsapce of $\Hom(N, H_1^\ast \otimes
V^\ast)$ defined by this condition.   Then $u_1 \mapsto \omega$ defines an
isomorphism
\[
\Hom(N, H_1^\ast \otimes V^\ast)^\prime \cong R(\bar{\omega})
\]
because $u_1 \mapsto \omega$ is injective since $\varphi$ and $\bar{u}$ are
surjective and because any $\omega \in R(\bar{\omega})$ arises in this way.
On the other hand, the factorization of $\tilde{u}_1 : N \otimes \ko_\P
\twoheadrightarrow \kc_{\bar{\omega}} \xrightarrow{\gamma} H_1^\ast \otimes
\ko_\P (1)$ defines an isomorphism $u_1 \leftrightarrow \gamma$ between
$\Hom(N, H_1^\ast \otimes V^\ast)^\prime$ and $\Hom\bigl(\kc_{\bar{\omega}},
H^\ast_1 \otimes \ko_\P (1)\bigr)$.   

b) By a) any $\omega\in R(\bar{\omega})$ gives rise to a diagram
\[
\xymatrix{& & & 0\\
& & &\bar{H}^\ast\otimes\ko_\P(1)\ar[u]&\\
0 \ar[r] & \bar{H}\otimes \ko_\P(-1)\ar[r]^{\bar{\beta}} & N\otimes
\ko_\P\ar[ur]^{\bar{\alpha}}\ar[r]\ar[d]^{\alpha_1} &
\kc_{\bar{\omega}}\ar[u]\ar[dl]\ar[r] & 0\\
& & H_1^\ast\otimes \ko_\P(1) & \ke_{\bar{\omega}}\ar[l]_-\pi\ar[u]\\
& & & 0\ar[u]}
\]
with $\pi$ induced by $u_1$ or $\alpha_1$. It follows that
$\alpha=\alpha_1+\bar{\alpha}: N\otimes \ko_\P\to (H_1^\ast\oplus
\bar{H}^\ast)\otimes \ko_\P(1)$ is surjective if and only if $\pi$ is
surjective. This proves b), because the surjectivity of $\alpha$ (as right
part of the monad of $\omega$) is equivalent to the non--degeneracy of
$\omega$. In that case the induced sequence
\[
H_1\otimes \ko_\P(-1)\to\ke_{\bar{\omega}}\xrightarrow{\pi} H_1^\ast\otimes
\ko_\P(1)
\]
is also a monad for the bundle $\ke_\omega$.
\end{proof}

\begin{subsub}\label{locfib} {\bf Corollary}: If $M^0(n-1,r+2)$ is irreducible
  of the expected dimension $3(n-1)n-\binom{2n-4-r}{2}$, then also
  $M(n,r)'_\xi$ is irreducible of the expected dimension
\[
3n(n+1)-\binom{2n-r}{2},
\]
if it is not empty.

\begin{proof} The dimension of any component of $M(n,r)$ is $\geq$ then the
  expected dimension by \ref{expd}. Because the fibres of $\res_\xi$ have
  constant dimension by \ref{fibresxi}, the corollary follows from Lemma \ref{cricomp}.
\end{proof} 
\end{subsub}
\end{sub}
\vskip5mm

\begin{sub}\label{union} {\bf Lemma}: If $M(n,r)'_\xi\neq \emptyset$, 
then for any other $\eta\neq 0,$
\[
M(n,r)'_\xi\cap M(n,r)'_\eta\neq \emptyset.
\]
 
\begin{proof} For $\omega_0\in M(n,r)'_\xi$ we consider the set 
\[
U_0=\{\eta\in H^\ast\ |\ \Im(\omega_0)\cap (\eta\otimes V^\ast)=0,\
h^1\ke_\eta(1)=0\}
\]
where $\ke_\eta$ denotes the bundle obtained from $\omega^0$ by restriction to
the kernel of $\eta$. Because $(\ke_\eta)$ is a flat family on the open set of
$\eta$'s defined by $\Im(\omega_0)\cap (\eta\otimes V^\ast)=0$, the
semicontinuity theorem implies that also $U_0$ is open. Because $\xi\in U_0$,
there is a $\xi'\in U_0$ which is independent of $\xi$. Then 
$\omega_0\in M(n,r)'_\xi\cap M(n,r)'_{\xi'}$.

Let now $\eta\neq 0$ be arbitrary in $H^\ast$. There is a transformation $g\in
GL(H)$ such that $g^\ast\xi=\xi$ and $g^\ast \xi'=\eta$. Then
\[
(g^\ast\otimes \id_{V^\ast})\circ\omega_0\circ (g\otimes \id_V)\in
M(n,r)'_\xi\cap M(n,r)'_\eta.
\] 
\end{proof}
\end{sub}

\begin{sub}\label{rest} \rm {\bf The induction}: Suppose now that $M^0(n-1,
  r+2)$ is irreducible and that also $M(n,r)'_\xi\neq \emptyset$ for some
  $\xi\neq 0$. Then by Lemma \ref{union} the union
\[
\underset{\xi\neq 0}{\cup} M(n,r)'_\xi
\]
is an irreducible open subset of $M(n,r)$ of the expected dimension. If,
in addition, $M^0(n-1, r+2)$ is smooth (as transversal intersection), then by
Lemma \ref{vanind} also $H^2S^2\ke_\omega=0$ for any $\omega$ in the above
union.
Concerning rank $2$ instantons, we shall prove that
\[
M(n,2)=\underset{\xi\neq 0}{\cup} M(n,2)_\xi
\]
for $n\geq 3$. For $n=3,4$ we shall even prove that 
\[
M(n,2)=\underset{\xi\neq 0}{\cup} M(n,2)'_\xi.
\]
This will give a unified proof of the instanton conjecture for $n\leq 4$.

For $n=5$, however, we are at present not able to prove that for any
$\omega\in M(5,2)$ there exists a $\xi\neq 0$ such that
$H^1\ke_{\bar{\omega}}(1)=0$. But we shall prove the weaker result, that there
exists a $\xi\neq 0$ with $h^1\ke_{\bar{\omega}}(1)\leq 1$. This is already
sufficient to prove the conjecture for $n=5$ in the sequel.

For $n=6$ one might hope that $\cup M(6,2)'_\xi$ is the complement of the
subvariety of 'tHooft bundles (defined by $h^0\ke(1)\neq 0$). However, for
$n>6$ the present induction method doesn't seem to work anymore.
\end{sub}
\vskip1cm

\section{On Jumping Lines of $n$--Instantons}\label{section4}

In this section we are going to prove the following Proposition \ref{exxi}
which enables us to choose $\xi\in H^\ast_n$ for any $n$--instanton
$\ke_\omega, n\geq 3$, such that $\ke_\xi:=\ke_{\bar{\omega}}$ is again an
$(n-1,4)$--instanton, i.e.\ $H^0\ke_\xi=0$. For $n\geq 5$ the proposition
allows us to choose a second $\eta\in H^\ast_{n-1}$ such that also
$H^0(\ke_\xi)_\eta=0$

In order to prepare the proof we include the following lemmata on jumping
lines of arbitrary instantons.

\begin{sub}\label{2lin} {\bf Lemma}: Let $\ke$ be a stable rank $2$ vector
  bundle on $\P_2$ with Chern class $(c_1, c_2)=(0,n),\ n \geq 2$, let $L_1,
  L_2$ be distinct lines and let $a_1, a_2\geq 0$ be defined by
  $\ke_{L_\nu}\cong\ko_{L_\nu}(a_\nu)\oplus \ko_{L_\nu}(-a_\nu)$. Then
  \begin{enumerate}
  \item [(a)] $a_1+a_2\leq n$
  \item [(b)] If $a_1+a_2=n$, and $a_1 \geq 2$, $a_2 \geq 2$, then $\ke$ can
    be realized as an extension
\[
0\to\ko_{\P_2}(-1)\to\ke\to\ki_Z(1)\to 0
\]
with $Z\subset L_1 \cup L_2$, a $0$--dimensional subscheme of length $n+1$.
  \end{enumerate}

\begin{proof} (a) There is the natural exact sequence 
$$
0\to \ko_{L_1\cup L_2}\to \ko_{L_1}\oplus \ko_{L_2}\to \ko_{L_1\cap L_2}\to
0.\eqno(1)
$$
Tensoring this with $\ke$ and taking sections, this implies
\[
a_1+a_2\leq h^0\ke_{L_1\cup L_2}.
\]
On the other hand, the sequence
$$
0\to\ke(-2)\to \ke\to \ke_{L_1\cup L_2}\to 0\eqno(2)
$$
implies the exact sequence
\[
0\to H^0\ke\to H^0\ke_{L_1\cup L_2}\to H^1\ke(-2)\to H^1\ke.
\]
Because $\ke$ is stable, $h^0\ke=0$ and $h^2\ke(-2)=h^0\ke(-1)=0$, so that
$h^1\ke(-2)=n$ by the Riemann--Roch formula, and then $h^0\ke_{L_1\cup
  L_2}\leq n$.

(b) According to the proof of (a), if $a_1+a_2=n$, then
$a_1+a_2=h^0\ke_{L_1\cup L_2}=n$ and
\[
H^0\ke_{L_1}\oplus H^0\ke_{L_2}\to H^0\ke_{L_1\cap L_2}
\]
is surjective. It follows that then also
\[
H^0\ke_{L_1}(1)\oplus H^0\ke_{L_2}(1)\to H^0\ke_{L_1\cap L_2}(1)
\]
is surjective and from this, that 
\[
h^0\ke_{L_1\cup L_2}(1)=n+2.
\]
Using sequence (2) with $\ke(1)$, we obtain $h^0\ke(1)\geq 2$. Then there are
exact sequences
\[
0\to \ko_{\P_2}(-1)\to\ke\to \ki_Z(1)\to 0
\]
with $0$--dimensional subschemes $Z$ of $\P_2$ of length $n+1$. Tensoring (1)
with $\ko_Z$ and taking lengths, one obtains
\[
\begin{array}{ccc}
\length(Z\cap (L_1\cup L_2))& \geq & \length(Z\cap L_1) +
\length(Z\cap L_2)-1\\
 & = & (a_1+1)+(a_2+1)-1\\
 & = & n+1 \\
 & = & \length(Z).
\end{array}
\]
Then $Z\subset L_1\cup L_2$ as schemes.
\end{proof}
\end{sub}
\vskip5mm

\begin{sub}\label{jlpl} {\bf Lemma}: Let $\ke$ be a semistable rank $2$ vector
  bundle on $\P_2$ with Chern class $(c_1, c_2)=(0,n),\ n\geq 2$. Then
  \begin{enumerate}
  \item [(a)] For odd $n\geq 3$ or even $n\geq 6$, $\ke$ has only finitely
    many jumping lines of order $\geq n/2$.
\item [(b)] For $n=2$ or $4$, $\ke$ has at most one jumping line of order $>n/2$.
  \end{enumerate}

\begin{proof} Note, that semistability in this case means that $H^0\ke(-1)=0$,
  see \cite{OSS}. Let firstly $\ke$ be stable with $H^0\ke=0$. Then
  \ref{2lin}, (a), implies that $\ke$ has at most one jumping line of order
  $>n/2$. Suppose $n\geq 4$ is even. According to \ref{2lin}, (b), if $\ke$
  has two jumping lines $L_1, L_2$ of order $\geq n/2$, then there is an
  extension
\[
0\to \ko_{\P_2}(-1)\to\ke\to\ki_Z(1)\to 0
\]
with $Z\subset L_1\cup L_2$ of length $n+1$. If $L$ is any other line, then
\[
\text{length } (Z\cap L)\leq \text{ length } (L\cap(L_1\cup L_2))=2,
\]
and it follows that $\ke_L\cong \ko_L(a)\oplus\ko_L(-a)$ with $0\leq a\leq
1$.

Now assume that $\ke$ is properly semistable. In this case we have an exact
sequence
\[
0\to \ko_{\P_2}\to\ke\to\ki_Z\to 0
\]
with a $0$--dimensional scheme $Z$ of length $n$. If $L$ is a jumping line of
$\ke$ order $\geq n/2$, then also length $(Z\cap L)\geq n/2$ by the exact
sequence tensored with $\ko_L$. If $\ke$ should have infinitely many jumping
lines of order $\geq n/2$, then there is a point $x\in Z$ such that infinitely
many such jumping lines meet $Z$ exactly in $x$. Then length
$(\ko_{Z,x}\otimes \ko_{L,x})\geq n/2$ for infinitely many lines. It follows
that the germ $(Z,x)$ is defined by two equations $f,g\in\ko_{\P_2,x}$ ($Z$ is
a locally complete intersection) with $\mult(f),\, \mult(g)\geq n/2$. Then
\[
n=\length(Z) \geq \length(\ko_{Z,x})\geq \mult(f)\mult(g).
\]
This is not possible for odd $n\geq 3$ or even $n\geq 6$. For $n=2$ or $4$
the same kind of argument as in the last part of the proof of \ref{2lin}, (b),
shows that $\ke$ has at most one jumping line of order $>n/2$.
\end{proof}
\end{sub}
\vskip5mm

\begin{sub}\label{jlsp} {\bf Proposition:} Let $\ke$ be an n-instanton on 
$\P_3$,\, $n\geq 2$. Then the set of jumping lines of $\ke$ of order $> n/2$
has dimension $\leq 1$.
\end{sub}

\begin{proof} Let $\G$ denote the Grassmannian of lines in $\P_3$ and consider
the incidence diagram 
$$\G \xleftarrow{p} \F \xrightarrow{q}\P_3^\ast $$
of lines in planes. Let $\Sigma$ be the set of jumping lines of order $>n/2$.
According to \ref{jlpl} the projection $q|p^{-1}(\Sigma)$ has finite fibres. 
Suppose $\dim\Sigma\geq 2$. Then $q(p^{-1}(\Sigma))=\P_3^\ast$. Let then $P_0$
be a stable plane for $\ke$. It must contain a jumping line $L_0$ of order 
$>n/2$. For a general plane $P$ containing $L_0$,\, $\ke_P$ is stable and, 
according to \ref{2lin}, (a), $L_0$ is the only jumping line of order 
$>n/2$ contained in $P$. On the other hand, for any  $P$ containing $L_0$ the
bundle $\ke_P$ is semistable. Hence, by \ref{jlpl}, $P$ contains at most finitely
many jumping lines of order $>n/2$. Consequently, there are only finitely many
jumping lines $L$ of order $>n/2$ meeting $L_0$. This means that 
$T_{L_0}(\G)\cap\Sigma$ is a finite set, where $T_{L_0}(\G)$ denotes the geometric
tangent hyperplane to $\G$ at $L_0$ in $\P_5$, whose intersection with $\G$ is
the cone of lines meeting $L_0$. This contradicts $\dim\Sigma\geq 2$.
\end{proof}
\vskip5mm

\begin{sub}\label{rmstr}{\bf Remark:}\rm\, One can show, using the method of 
  R.~Strano and M.~Green, as in \cite{C2}, that a non special 'tHooft
  n-instanton $\ke$,\, $n \geq 5$, satisfies $h^0\ke_P(1)\leq 1$ for a general
  plane $P$.  Then, using the above arguments, one deduces that the set of
  jumping lines of $\ke$ of order $\geq n/2$ is at most 1-dimensional. But
  this improves \ref{jlpl} only for even $n\geq 6$.
\end{sub}

From the above statements on jumping lines we can now derive the following 
proposition, which is the key of the induction process of this paper.
\vskip5mm

\begin{sub}\label{exxi} {\bf Proposition:} a) Let $\ke$ be any $n$--instanton. 
If $n\geq 3$, then
for a general $\xi\in H^1\ke(-1)$ the multiplication map 
$\xi\otimes  V^\ast\to H^1\ke$ is injective.

b) If $n\geq 5$, then for a general $2$--dimensional subspace $U\subset
H^1\ke(-1)$ the multiplication map $U\otimes V^\ast\to H^1\ke$ is
injective.
\end{sub}

Note, that the condition in a) means that $\Im(\omega)\cap (\xi\otimes
V^\ast)=0$, see \ref{rkl}, because the multiplication map is isomorphic to 
$H\otimes V^\ast\to Q$.

\begin{proof} a) The right part $N\otimes \ko\to H^\ast\otimes \ko(1)$ of the
  monad of $\ke$ corresponds to the exact sequence $0\to N\to H^\ast\otimes
  V^\ast\xrightarrow{\mu} Q\to 0$, where $\mu$ is isomorphic to the
  multiplication map. If $P\subset\P$ is a plane with equation $z\in V^\ast$,
  we obtain the exact sequence
\[
0\to H^0\ke_P\to N\to H^\ast\otimes V^\ast/\langle z\rangle 
\]
and, therefore,
$$
H^0\ke_P\cong N\cap (H^\ast\otimes z).\eqno (1)
$$

On the other hand, if $h^0\ke_P\neq 0$, we are given an exact sequence
\[
0\to \ko_P\to\ke_P\to\ki_Z\to  0
\]
where $Z$ is a $0$--dimensional subscheme of  $P$ of length $n$, because
$h^0\ke_P(-1)=0$. It follows that
$$
0\leq h^0\ke_P\leq 1. \eqno (2)
$$
We let now $\P(N)_1$ denote the set of decomposable classes $\langle
\xi\otimes z\rangle \in\P(H^\ast\otimes V^\ast)$ which are contained in
$P(N)$, so that
\[
\P(N)_1=\P(N)\cap S_1
\]
where $S_1$ is the image of the Segre embedding $\P H^\ast\times \P
V^\ast\subset\P(H^\ast\otimes V^\ast) $. There are the two projections
\[
\P H^\ast\xleftarrow{p_1}\P(N)_1\xrightarrow{q_1}\P V^\ast.
\]
The isomorphism (1) implies that for any $z\neq 0$ we have
\[
\P(H^0\ke_P)\cong q_1^{-1}(\langle z\rangle ),
\]
where $P$ is the plane with equation $z$. Therefore, the image of $q_1$ is
contained in the subvariety $JP(\ke)\subset \P V^\ast$ of unstable planes of
$\ke$. 
\vskip2mm

Now (2) implies that $\P(N)_1$ is isomorphic to $JP(\ke)$. Because $\dim
JP(\ke)\leq 2$, see \cite{Barth1}, it follows that $p_1(\P(N)_1)$ has dimensions $\leq 2$. This
proves a) in case $n\geq 4$. If $n=3$ and $\ke$ is not special 'tHooft, then
$\dim JP(\ke)\leq 1$ by \ref{fures}, (4), so that a) is true in that case,
too. But in case $n=3$ and $\ke$ is special 'tHooft, the claim follows from
remark \ref{remsub} below.

b) For the proof of part b) we consider the intersection
\[
\P(N)_2=\P(N)\cap S_2
\]
with the secant variety of $S_1$ in $\P(H^\ast\otimes V^\ast)$, such that
$\P(N)_2$ consists of all elements of type $\langle\xi_1\otimes z_1 +
\xi_2\otimes z_2\rangle$ which are contained in $\P(N)$. Then we have the two
projections
\[
G(2, H^\ast)\xleftarrow{p_2}\P(N)_2\smallsetminus \P(N)_1\xrightarrow{q_2}
G(2, V^\ast).
\]
We have to show that a general $2$--dimensional subspace $K\subset H^\ast$
satisfies $N\cap (K\otimes V^\ast)=0$.

Now $N\cap (K\otimes V^\ast)\neq 0$ if and  only if either $\P(K)\cap
p_1(\P(N)_1)\neq \emptyset$ or $K\in\Im(p_2)$.

Because $\dim p_1(\P(N))\leq 2$ by part a), it follows that the set of $K$ with
$\P(K)\cap p_1(\P(N)_1)\neq \emptyset$ has dimension $\leq 2+(n-2)=n$.  By
assumption, $\dim G(2, H^\ast)=2(n-2)>n+1$ for $n\geq 6$ and $=n+1$ for
$n=5$. Therefore, it remains to prove that $\dim \Im(p_2)\leq n+1\ \text{ for
  } n\geq 6$ and $\dim \Im(p_2)\leq 5\ \text{ for  } n=5$.

In order to derive the estimates, we consider the fibres of $q_2$. Let
$W\subset V^\ast$ be any $2$--dimensional subspace with the dual line
$L\subset \P V$. Then
\[
N\cap (H^\ast\otimes W)\cong H^0\ke_L
\]
as can easily be derived from the monad description of $\ke$. On the other
hand,
\[
q_2^{-1}(W)=\P(N\cap H^\ast\otimes W)\smallsetminus\P(N)_1
\]
by definition of $\P(N)_2$. Therefore, 
$$
\dim q_2^{-1}(W)=h^0\ke_L-1.\eqno(3)
$$
Let $G(2, V^\ast)\supset J_1\supset J_2\supset\cdots\supset J_n$ be the
filtration by the sets $J_k$ of jumping lines $L$ of order $h^0\ke_L-1\geq
k$. It follows that
\[
\dim q_2^{-1}(J_k\smallsetminus J_{k+1})\leq 3+k\leq  n
\]
for $k\leq n-3$. For $k=n-2$
\[
\dim q_2^{-1}(J_{n-2}\smallsetminus J_{n-1})\leq 1+(n-2)=n-1
\]
because of Lemma \ref{jlsp} since $n-2>n/2$, and finally
\[
\dim q_2^{-1}(J_{n-1}\smallsetminus J_n)\leq n\quad ,\quad \dim
q_2^{-1}(J_n)\leq n+1.
\]
Totally we have $\dim \bigl(\P(N)_2\smallsetminus \P(N)_1\bigr)\leq n+1$,
which is sufficient for $n\geq 6$. In case $n=5$, $\dim
\bigl(\P(N)_2\smallsetminus \P(N)_1\bigr)=6$ is only possible if $\dim J_5=1$.
In that case $\ke$ is a special 'tHooft bundle, and it follows by the direct
argument in remark \ref{remsub}, that there are $2$--dimensional subspaces
$K\subset H^\ast$ with $N\cap(K\otimes V^\ast)=0$. This proves b).
\end{proof}

\vskip5mm

\begin{sub}\label{remsub}{\bf Remark:}\rm\, Let $\ke$ be an $(n,2)$--instanton on
  $\P_3$ with $n\geq 2m+2\geq 3$. One can prove, along the same lines, that
  then, for a general $m$--dimensional subspace $K\subset H^1\ke(-1)$ the
  multiplication map $K\otimes V^\ast\to H^1\ke$ is injective. If $\ke$ is a
  special 'tHooft bundle, this can be verified by the following direct
  argument. Let $z_0, \ldots, z_3$ be a basis of $V^\ast$. Then bases of $N$
  and $H^\ast$ can be chosen such that the operator $N\to H^\ast\otimes
  V^\ast$, i.e.\ the right part of the monad of $\ke$, can be represented by
  the matrix
\[
\left(
\begin{array}{llllllll}
z_0 & z_1 & z_2 & z_3 & \\
    &     & z_0 & z_1 & z_2 & z_3 & \\
& &  &\ddots &  & &\ddots\\
& & & & z_0 & z_1 & z_2 & z_3
\end{array}
\right).
\]
If $n=2m+1$ or $2m+2$ and $e_1, \ldots, e_n$ is the corresponding basis of
$H^\ast$, then\\
$K=\span(e_2, e_4, \ldots, e_{2m})$
satisfies $N\cap (K\otimes V^\ast)=0$.
\end{sub}

For later use, we consider the space $\Pi (\xi) \subset G(2, H^\ast)$ of all
2--dimensional subspaces $K \subset H^\ast$ which contain $\xi$.   We have:
\vskip5mm

\begin{sub}\label{indec} {\bf Lemma:} Let $n = 5$ and $\omega \in M(5,2)$ such that
  $\ke_\omega$ is not special 'tHooft.  Then for a general $\xi \in H^\ast$
  there are closed subsets $T_1 \subset T_2 \subset \Pi(\xi)$, $\dim T_1 \le
  1$ and $T_2$ a surface, such that

  \begin{enumerate}
  \item[(i)] $N \cap (K \otimes V^\ast) = 0$ for $K \in \Pi(\xi)
    \smallsetminus T_2$;

  \item[(ii)] $\dim N \cap (K \otimes V^\ast) = 1$ and $N \cap (K \otimes
    V^\ast)$ contains no non--zero decomposable vector of $H^\ast \otimes
    V^\ast$ for $K \in T_2 \smallsetminus T_1$.  
  \end{enumerate}
\end{sub}
\vskip5mm

\begin{proof}
  We use the previous notation and consider the morphisms
\[
\P(N)_1 \xrightarrow{p_1} \P(H^\ast) \text{ and } \P(N)_2 \smallsetminus
\P(N)_1 \xrightarrow{p_2} G(2, H^\ast) =: G.
\]

a)\; Let $Y_1 \subset G$ denote the subset of those $K \in G$ for which $N
\cap (K \otimes V^\ast)$ contains a non--zero decomposable vector or,
equivalently, $\P(K) \cap p_1\bigl(\P(N)_1\bigr) \not= \emptyset$.   Because
$\dim p_1 \bigl(\P(N)_1\bigr) \le 1$, the set of lines $\P(K)$ in $\P(H^\ast)$
with this condition is closed and of dimension $\le 4$.

b)\; Let $Y_2 \subset G$ be the set of all $K \subset G$ for which $N \cap (K
\otimes V^\ast) \not= 0$.   If $\kk \subset H^\ast \otimes \ko_G$ denotes the
universal subbundle, we consider the homomorphism
\[
\kk \otimes V^\ast \xrightarrow{\phi} Q \otimes \ko_G,
\]
obtained as the composition of $\kk \otimes V^\ast \subset 
H^\ast\otimes V^\ast\otimes\ko_G$ 
and the multiplication map $H^\ast \otimes V^\ast \to Q$ with kernel
$N$.  Then $Y_2$ is the determinantal subvariety of $\phi$.   Because, for a
general $K$ we have $N \cap (K\otimes V^\ast) = 0$ by Proposition \ref{exxi},
$Y_2$ is a hypersurface in $G$, $\dim Y_2 = 5$.   By definition of $Y_2$ we
have $\Im(p_2) \subset Y_2$, and, furthermore,
\[
Y_2 = Y_1 \cup \Im(p_2).
\]

c)\; For $K \in Y_2 \smallsetminus Y_1$ the fibre $p_2^{-1}(K)$ consists of
vectors $\xi_1 \otimes z_1 + \xi_2 \otimes z_2 \in N$, where $\xi_1, \xi_2$ is
a basis of $K$ and $z_1, z_2 \in V^\ast$, i.e.
\[
p_2^{-1}(K) \cong \P\bigl(N \cap (K \otimes V^\ast)\bigr).
\]
Now, by the proof of Proposition \ref{exxi}, $p_2\bigl(\P(N)_2 \smallsetminus
\P(N)_1\bigr)$ has dimension $\le 5$.   It follows that, for a general $K$ in
each component of $Y_2 \smallsetminus Y_1$, the fibre $p_2^{-1}(K)$ is a point,
or, equivalently, $N \cap (K \otimes V^\ast)$ is 1--dimensional.   Then the
subvariety $Y_2^\prime \subset Y_2$, defined by $\dim N \cap (K \otimes
V^\ast) \ge 2$, has dimensional $\le 4$.   Let $Z_1 = Y_1 \cup Y_2^\prime$.
Then $\dim Z_1 \le 4$ and for $K \in Y_2 \smallsetminus Z_1$ the intersection
$N \cap (K \otimes V^\ast)$ is 1--dimensional and contains no non--zero
decomposable vector.

d)\; Because $\dim Z_1 \le 4$, the general 3--space $\Pi(\xi)$ meets $Z_1$
at most in dimension 1.   Let, now, $T_1 = \Pi (\xi) \cap Z_1$ and $T_2 =
\Pi(\xi) \cap Y_2$.   We may assume that $\Pi(\xi) \not\subset Y_2$ by
Proposition \ref{exxi}, b), so that $T_2$ is a surface in $\Pi(\xi)$.   Then
$T_1 \subset T_2$ satisfy the properties of the lemma. 
\end{proof}
\vskip5mm

\section{The varieties $M(n,r)$ for $n\leq 4$}\label{section5}

We first note that $M(n,2n)$ is an open subset of the affine space
$S^2H^\ast\otimes \Lambda^2 V^\ast$. In fact, it is the complement of the
hypersurface $M_{4n-2}$. In this case, if $\omega\in M(n,2n)$, the bundle
$\ke_\omega$ is the cokernel in
\[
0\to H \otimes \ko_\P(-1)\to H^\ast\otimes \Omega^1(1)\to\ke_\omega\to 0,
\]
see \ref{beidia} with $Q=0$. Then $H^1\ke_\omega(i)=0$ for $i\geq 0$ and we have
\[
M^0(n,2n)=M(n,2n).
\]
\vskip5mm

\begin{sub}\label{m6} {\bf Proposition:} For $n=2,\, M_6=M_{4n-2}(H)$ is an irreducible
  hypersurface in $S^2H^\ast\otimes \Lambda^2 V^\ast$.

\begin{proof} Note that $M_6\smallsetminus (M_4\cup\Delta)$ had been shown in
  \ref{hcase} to be a smooth hypersurface, and $M_6$ is a homogeneous
  hypersurface in $S^2H^\ast_2\otimes \Lambda^2 V^\ast$.

In order to prove that it is irreducible, it suffices to prove that it is
non--singular in codimension 1. Now $\Delta=\Delta(H_2)$ is closed and irreducible in $S^2H^\ast_2\otimes \Lambda^2V^\ast$ of codimension 2. We have
$M_4\subset\Delta$ and hence $M_6\smallsetminus \Delta$ is smooth. It is
therefore sufficient to find an $\omega\in\Delta\smallsetminus M_4$ which is a
smooth point of $M_6$. For that choose an isomorphism $H_2\cong k^2$ and let
$H_2\xrightarrow{\omega} H_2^\ast\otimes \Lambda^2 V^\ast$ be presented by the
matrix
\[
\left(
  \begin{array}{lc}
x_2\wedge x_3 & x_2\wedge x_4\\
x_2\wedge x_4 & x_1\wedge x_2+x_3\wedge x_4
  \end{array}\right)
\]
where $x_1, \ldots, x_4$ is a basis of $V^\ast$. One can easily check that
this $\omega$ is degenerate at the point $x_2=x_3=x_4=0$ and that
$\rk(\omega)=6$. Moreover, $\omega$ is a smooth point of
$M_6$, using the argument of Katsylo--Ottaviani as in the proof of Proposition
\ref{hcase}.
\end{proof}
\end{sub}
\vskip5mm

\begin{sub} {\bf Corollary 1:} {\rm (Hartshorne)}
$MI(2)$ is smooth and irreducible of the expected dimension $13$.

\begin{proof}$M(2,2)=M_6\smallsetminus (M_4 \cup\Delta)$ is a smooth
  transversal intersection and irreducible. By \ref{conj2} the results follows.
\end{proof}
\end{sub}
\vskip5mm

\begin{sub}\label{ellstr} {\bf Corollary 2:} {\rm (Ellingsrud - Str{\o}mme)}
$MI(3)$ is smooth and irreducible of the expected dimension 21.

\begin{proof} As noted above $M^0(2,4)=M(2,4)$, and this is irreducible. By
  Proposition \ref{exxi} $M(3,2)$ is the union of the open sets $M(3,2)_\xi$
  which are equal to $M(3,2)'_\xi$ and which are transversal intersections
  and irreducible. The result follows now as in the previous case.
\end{proof}
\end{sub}

In order to treat the case $MI(4)$, we prove the following 

\vskip5mm

\begin{sub}\label{m10}{\bf Lemma:} For $n=3,\, M_{10}=M_{4n-2}(H)$ is an irreducible
  hypersurface in $S^2H^\ast\otimes \Lambda^2 V^\ast$.
\end{sub}
\begin{proof} By \ref{hcase} $M_{10}$ is a hypersurface in $S^2H^\ast_3\otimes
  \Lambda^2 V^\ast$ and $M(3,4)=M_{10}\smallsetminus (\Delta\cup M_8)$ is
  smooth. Because for $n=3$ we have $M_6\subset\Delta,\ 
  M(3,2)=M_8\smallsetminus \Delta$ and hence $\Delta\cup M_8=\Delta\cup
  M(3,2)$. Now $\Delta$ is irreducible in $S^2H^\ast_3\otimes \Lambda^2
  V^\ast$ of codimension $4$, see remark in \ref{rkstr}, and $M(3,2)$ has
  codimension $6$ in $S^2H^\ast_3\otimes \Lambda^2 V^\ast$ by corollary
  \ref{ellstr}. It follows that $M_{10}$ is smooth in codimension $1$ and so
  $M_{10}$ is irreducible.
\end{proof}
\vskip5mm

\begin{sub}\label{three4}{\bf Corollary:}\ $M(3,4)$ is irreducible and a (smooth) transversal
  intersection of expected dimension 45.
\end{sub}
\vskip5mm

\begin{sub}{\bf Proposition:} {\rm (Barth, LePotier)}
$MI(4)$ is smooth and irreducible of the expected dimension $29$.

\begin{proof} 1) Recall that the open subsets $M(4,2)'_\xi\subset M(4,2)$ are
  defined by
\[
M(4,2)'_\xi=\{\omega\in M(4,2)\ |\ \bar{\omega}\in M(3,4)\text{ and }
H^1\ke_{\bar{\omega}}(1)=0\},
 \]
where $\bar{\omega}=\res_\xi(\omega)$. In this proof we write
$\kf=\ke_{\bar{\omega}}$, which depends on the choice of $\xi$. By
\ref{indproc} $M(4,2)'_\xi$ is irreducible. It is also smooth: by \ref{crism}
$H^2 S^2\kf=0$ for $\bar{\omega}\in M(3,4)$ because $M(3,4)$ is a smooth
transversal intersection of the expected (co)dimension. It follows from Lemma
\ref{vanind} that also $H^2S^2\ke_\omega=0$ for $\omega\in M(4,2)'_\xi$, hence
again by \ref{crism} $M(4,2)'_\xi$ is smooth at any of its points and of
expected dimension.

2) It suffices now to prove that
\[
M(4,2)=\underset{\xi\neq 0}{\cup} M(4,2)'_\xi.
\]
Let $\omega\in M(4,2)$. By \ref{vanpl} there is $\xi\in H^1\ke_\omega(-1)$
such that $H^1\kf_P(1)=0$ for any plane $P$ in $\P_3$, provided $\ke_\omega$ is
not special 'tHooft. In case $n=4$ we obtain $h^1\kf=2$. Then the operator
$H^1\kf\otimes V^\ast  \to H^1\kf(1)$ can be presented by a matrix
\[
A=\left(
  \begin{array}{lcl}
v_{11} & \ldots & v_{1m}\\
v_{21} & \ldots & v_{2m}
  \end{array}
\right)
\]
with $m=h^1\kf(1)$ if $h^1 \kf(1)\neq 0$. Let now $P$ be a plane with equation
$z=0$ and such that $z(v_{11})=z(v_{21})=0$. Then
\[
H^1\kf\xrightarrow{A(z)} H^1\kf(1)
\]
cannot be surjective and it would follow that $H^1\kf_P(1)=\coker A(z) \neq
0$, a contradiction. This proves that $\omega\in M(4,2)'_\xi$ for the chosen
$\xi$. If, however, $\ke_\omega$ is special 'tHooft, the pairing
$H^1\ke_\omega(-1)\otimes V^\ast\to H^1\ke_\omega$ can be presented by a
$4\times 6$ matrix
\[
\left(
\begin{array}{lll}
v_1 v_2 & &\\
v_3 v_4 & v_1 v_2 & \\
        & v_3 v_4 & v_1 v_2\\
        &         & v_3 v_4
\end{array}
\right)\raisebox{-4ex}{,}
\]
where $v_1, \ldots, v_4$ is a basis of $V$, see \cite{BoeTrm}. Then,
choosing $\xi=(0,1,0,0)\in k^4\cong H^1\ke_\omega(-1)$, the resulting
homomorphism $H^1\kf(-1)\otimes V^\ast\to H^1\kf$ is presented by the matrix
\[
A=\left(
  \begin{array}{cc}
0 & 0\\
v_1 & v_2\\
v_3 & v_4
  \end{array}\right)\raisebox{-4ex}{.}
\]
It follows that $H^1\kf(1)=0$ because of the exact sequence $H^1
\kf(-1)\otimes \Lambda^2V^\ast\to H^1\kf\otimes V^\ast\to H^1\kf(1)\to 0$.
\end{proof}
\end{sub}
\vskip5mm

\begin{sub}\label{rm44}{\bf Remark}: \rm Concerning $MI(5)$ and $M(4,4)$, it
follows from \ref{indproc} that the open set
\[
\underset{\xi\neq 0}{\cup} M(4,4)_\xi
\]
is irreducible and a (smooth) transversal intersection of expected dimension,
because $M^0(3,6)=M(3,6)$ is smooth of expected dimension. This will suffice
to prove that $MI(5)$ is smooth and irreducible of dimension $37$. Using the
method of Katsylo--Ottaviani as in the proof of \ref{hcase}, one can show that
$M(4,4)$ is smooth of the expected dimension. It is an open question whether
it is also irreducible.
\end{sub}
\vskip1cm

\section{A technical result about $5$--instantons}\label{section6}

\begin{sub}\label{boundh1}{\bf Proposition}: Let $\ke_\omega$ be a 
$5$--instanton on $\P_3$. Then for a
  general $\xi\in  H^1\ke_\omega(-1)$ the associated rank--$4$ bundle
  $\kf=\ke_{\bar{\omega}}$ satisfies $h^1\kf(1)\leq 1$.
\end{sub}\vskip5mm

For the proof we need the following lemmata on the vanishing of $H^1\kf_L(1)$
and $H^1\kf_P(1)$ for lines and planes.
\vskip5mm

\begin{sub}\label{vanli} {\bf Lemma}: Let $\ke_\omega$ be a $n$--instanton,
  $n=4$ or $5$,  and
  assume that $\ke_\omega$ is not a special 'tHooft bundle
  $(h^0\ke_\omega(1)<2)$. Then, for a general $\xi\in H^1\ke_\omega(-1)$, the
  bundle $\kf=\ke_{\bar{\omega}}$ has the property that $H^1\kf_L(1)=0$ for
  any line $L$ in $\P_3$. 
\end{sub}

\begin{proof} Let $\ke=\ke_\omega$. For any plane $P$ in $\P_3$ we have
  $H^1\ke(-1)\xrightarrow{\approx} H^1\ke_P(-1)$. Then for any line
  $L\subset P$ the restriction $H^1\ke(-1)\to H^1\ke_L(-1)$ is surjective
  because $H^1\ke_P(-1)\to H^1\ke_L(-1)$ is surjective, since
  $H^1\ke_P(-2)=0$. Because $\ke_L$ is the cohomology of the monad
\[
0\to \ko_L(-1)\to\kf_L\to \ko_L(1)\to 0,
\]
we obtain  
\[
H^1\kf_L(1)\cong H^1\ke_L(1)/\xi_L\cdot H^0\ko_L(2),
\]
where $\xi_L$ denotes the restriction of $\xi$ in $H^1\ke_L(-1)$. On the other
hand, $\ke_L\cong\ko_L(a)\oplus\ko_L(-a),\ 0\leq a\leq n$. If $0\leq a\leq 2$,
then $H^1\ke_L(1)=0$. If $a=3, H^1\ke_L(-1)\cong H^1\ko_L(-4)$ and
$H^1\ke_L(1)\cong H^1\ko_L(-2)$. Then for $\eta\in H^1\ko_L(-4)$ we have
$\eta\cdot H^0\ko_L(2)\neq H^1\ko_L(-2)$ if and only if
$\eta=0$. Consequently, for the vanishing of $H^1\kf_L(1)$  we have only to
assume that $\xi_L\neq 0$ or equivalently that
\[
\xi\not\in\Ker(H^1\ke(-1)\to H^1\ke_L(-1)),
\]
which is $1$-- or $2$--dimensional. Since by \ref{jlsp} the set of jumping lines of
$\ke$ of order $\geq 3$ is at most $1$--dimensional, $\xi$ has to avoid a
subvariety of dimension $\leq 2$ or $3$. If $a=4$, the elements $\xi$ or $\eta\in
H^1\ko_L(-5)$ should avoid the condition $\eta H^0\ko_L(2)\neq
H^1\ko_L(-3)$. The set of these $\eta$ is the affine cone over the rational
normal curve in $\P H^1\ko_L(-5)$. Namely, if $s,t$ are homogeneous coordinates
on $L$, we have
\[
\eta=\sum\limits_{\nu=1}^4\frac{a_\nu}{s^{5-\nu}t^\nu}\quad\text{ in }\quad
H^1\ko_L(-5)
\]
and
\[
s^2\eta=\frac{a_1}{s^2t}+\frac{a_2}{st^2},\quad
st\eta=\frac{a_2}{s^2t}+\frac{a_3}{st^3},\quad
t^2\eta=\frac{a_3}{s^2t}+\frac{a_4}{st^2}.
\]
Then the condition $\eta H^0\ko_L(2)\neq H^1\ko_L(-3)$ becomes 
\[
\rk
\left(
\begin{array}{lll}
a_1 & a_2 & a_3\\
a_2 & a_3 & a_4
\end{array}
\right)\leq 1.
\]
So $\xi$ has to avoid another $3$--dimensional subvariety in
$H^1\ke(-1)$. Finally, if $a=5$, in case $n=5$ only, the set of $\eta\in H^1\ko_L(-6)$ with $\eta
H^0\ko_L(2)\neq H^1\ko_L(-4)$ is a cubic hypersurface in $H^1\ko_L(-6)$, as
can be seen by a similar argument. So $\xi$ has to avoid this hypersurface
in $H^1\ke(-1)\cong H^1\ke_L(-1)$ in case $a=5$. But $\ke$ has only finitely
many jumping lines of order $5$ because it is not special 'tHooft, see also
\cite{Skiti}, \cite{Rao}. Totally
$\xi$ has to avoid two $3$--dimensional and finitely many $4$--dimensional
subvarieties in $H^1\ke(-1)$, in case $n=5$,  in order to satisfy the condition of the lemma. 
\end{proof}
\vskip5mm

\begin{sub}\label{vanpl} {\bf Lemma}: Let $\ke=\ke_\omega$ be as in the
  previous lemma. Then, for a general $\xi\in H^1\ke_\omega(-1)$, the bundle
  $\kf=\ke_{\bar{\omega}}$ has the property that $H^1\kf_P(1)=0$ for any plane 
$P$ in $\P_3$.
\end{sub}
\begin{proof} By \ref{vanli} we may assume that $H^1\kf_L(1)=0$ for any line
  $L$ if $\xi$ is in a fixed open set of $H^1\ke(-1)$. We let $\xi_P$ denote
  the element corresponding to $\xi$ under the isomorphism $H^1\ke(-1)\cong
  H^1\ke_P(-1)$ for a plane $P$. We are going to show that $H^1\kf_P(1)=0$ for
  any plane if $\xi$ avoids some additional subvarieties of $H^1\ke(-1)$. These
  subvarieties will be estimated in dimension in the following cases.

{\bf case 1}: $P$ is a stable plane of $\ke$. In this case there is no
condition on $\xi$ because we show that then already $H^1\kf_P(1)=0$.

Proof of case 1: We prove first that $h^1\kf_P\leq 2$. As in the previous
proof we have the exact sequence
\[
\xi_P\otimes H^0\ko_P(1)\xrightarrow{m(\xi)} H^1\ke_P\to H^1\kf_P\to 0.
\]
Because $h^0\ke_P=0$, we have $h^1\ke_P=n-2$. So, if $n=4$, then
$h^1\kf_P\leq 2$. If $n=5$, the homomorphism $m(\xi)$   is zero if
and only if $\xi_P$ is mapped to zero under the map $H^1\ke_P(-1)\to W\otimes
H^1\ke_P$, where $W^\ast=H^0\ko_P(1)$, and the kernel of this map is the image
of $\Hom(\Omega^1_P(1), \ke_P)$ in $\Ext^1(\ko_P(1), \ke_P)\cong
H^1\ke_P(-1)$. Therefore, $\xi_P$ is induced by a non--zero homomorphism
$\Omega^1_P(1)\xrightarrow{\varphi} \ke_P$. Because both bundles are stable,
$\varphi$ is generically and then globally injective. Then
$\kc=\coker(\varphi)$ has Hilbert polynomial $\chi\kc(m)=m-3$ and is
Cohen--Macaulay. It follows that $\kc=\ko_L(-4)$ for a line $L\subset P$. If
follows that $L$ is a jumping line of order $4$ and that the restricted
homomorphism $\Omega^1_P\otimes \ko_L\to\ke_L(-1)$ factors through
$\ko_L(3)$. Then the diagram 
\[
\xymatrix{H^1\Omega^1_P\ar[d]\ar[r] & H^1\ke_P(-1)\ar[d]\\
H^1(\Omega^1_P\otimes \ko_L)\ar[r] & H^1\ke_L(-1)
}
\]
implies that the image $\xi_L$ of $\xi$ in $H^1\ke_L(-1)$ is zero, and $H^1\kf_L(1)\cong
H^1\ke_L(1)/\xi_LH^0\ko_L(2)\cong H^1\ko_L(-3)\neq 0$
contradicting Lemma \ref{vanli}. This proves that $h^1\kf_P\leq 2$ if $\xi$ is
general.

Let now $L\subset P$ be any line with equation $z$. By the assumption
on $\xi$ the multiplication map $H^1\kf_P\xrightarrow{z} H^1\kf_P(1)$ is
surjective. Applying the bilinear map lemma of H.~Hopf to 
\[
H^1(\kf_P(1))^\ast\otimes H^0\ko_P(1)\to H^1(\kf_P)^\ast\ ,
\]
we deduce that 
\[
h^1\kf_P(1)\leq h^1\kf_P-h^0\ko_P(1)+1\leq 0.
\]
(when $K$ denotes the kernel of the map, the above surjectivity condition
implies that $\P H^1(\kf(1))^\ast\times \P H^0\ko_P(1)$ has an empty
intersection with $\P K$ in $\P(H^1(\kf(1))^\ast\otimes H^0\ko_P(1)),$ which
implies the estimate).

{\bf case 2}: $P$ is an unstable plane of $\ke$ but contains no jumping line
of order $5$. 

In this case we have an exact sequence
\[
0\to \ko_P\to\ke_P\to\ki_{Z, P}\to 0
\]
where $Z$ is a $0$--dimensional subscheme of $P$ of length $n$ with
$H^0\ki_{Z,P}(1)=0$, because $P$ contains no jumping line of order $n$. It 
follows that $h^0\ke_P=1,\ h^0\ke_P(1)=3$ and then $h^1\ke_P=n-1,\, 
h^1\ke_P(1)=n-3$. If
$h^1\kf_P\leq 2$, then, as in case 1, $h^1\kf_P(1)=0$. If $h^1\kf_P\geq 3$,
then the multiplication map $m(\xi)$, see case 1, has rank $n-4$. If $n=4$,
then $m(\xi)=0$, and we get a contradiction by the argument in case 1, which
leads to $h^1\kf_P\leq 2$. If $n=5$, then $m(\xi)$ has rank $1$ and $\xi_P$ is
annihilated by two linear forms $z_0, z_1\in H^0\ko_P(1)$.  
Let $x\in P$ be the point determined by $z_0, z_1$ 
and let $\ki_{\{x\}, P}$ be its
ideal sheaf. The standard resolution of this sheaf yields the exact sequence
\[
0\to H^0(\ke_P)^2\to H^0(\ki_{\{x\}, P}\otimes \ke_P(1))\to
H^1\ke_P(-1)\xrightarrow{z_0, z_1} H^1(\ke_P)^2.
\]
On the other hand, the defining sequence of $\ki_{Z,P}$ and $H^0\ko_P(1)\cong
H^0\ke_P(1)$ implies
\[
h^0(\ki_{\{x\}, P}\otimes \ke_P(1))=\left\{
  \begin{array}{lcl}
2 & \text{ if } & x\in P\smallsetminus Z\\
3 & \text{ if } & x\in Z
  \end{array}\right . .
\]
Therefore, $h^1\kf_P\geq 3$ can only occur if $\xi_P$ avoids at most five
$1$--dimensional vector subspaces of $H^1\ke_P(-1)$. Because $\ke$ has
only a $1$--dimensional variety of unstable planes, it follows that
$h^1\kf_P\leq 2$ and $h^1\kf_P(1)=0$ for any unstable plane, if $\xi$ avoids
an at most $2$--dimensional subvariety of $H^1\ke(-1)$.

{\bf case 3}: $P$ is an unstable plane of $\ke$ and contains a jumping line
$L$ of order $n$.

In this case $H^1\ke_P(-1)\xrightarrow{\approx} H^1\ke_L(-1)$ and
$h^1\ke_P(1)=n-2=h^1\ke_L(1)$, such that $H^1\ke_P(1)
\xrightarrow{\approx} H^1\ke_L(1)$. It follows that also
$H^1\kf_P(1)\xrightarrow{\approx} H^1\kf_L(1)=0$.
\end{proof}
\vskip5mm

{\bf Proof of proposition \ref{boundh1}:} Let us assume, firstly, that
$\ke=\ke_\omega$ is  not special 'tHooft. By \ref{exxi} we may assume that the
multiplication map $\xi\otimes V^\ast\to H^1\ke$ is injective, hence
$h^1\kf=4,\ \kf=\ke_{\bar{\omega}}$. By \ref{vanpl} we may assume that
$H^1\kf_P(1)=0$ for any plane. Therefore, the multiplication map $H^1\kf\to
H^1\kf(1)$ is surjective for any linear form. Now an application of the
bilinear lemma of H.~Hopf to
\[
H^1(\kf(1))^\ast\otimes V^\ast\to H^1(\kf)^\ast
\]
implies $h^1(\kf(1)\leq h^1\kf-4+1=1$. If $\ke$ is a special 'tHooft bundle,
it can even be shown that for a general $\xi\in H^1\ke(-1)$ we have
$H^1\kf(1)=0$. If $\ke$ is special 'tHooft, then the pairing $V^\ast\otimes
H^1\ke(-1)\to H^1\ke$ can be defined by a $5\times 8$ matrix
\[
\left(
  \begin{array}{llllllll}
v_1 & v_2 &     &     & &&&\\
v_3 & v_4 & v_1 & v_2 & & & &\\
    &     & v_3 & v_4 & v_1 & v_2 & &\\
    &     &     &     & v_3 & v_4 & v_1 & v_2\\
    &     &     &     &     &     & v_3 &v_4
  \end{array}\right)
\]
where $v_1, \ldots, v_4$ is a basis of $V$, after choosing suitable bases of
the cohomology spaces, see \cite{BoeTrm}. Then, choosing $\xi=(0,0,1,0,0)$ in
$ k^5\cong H^1\ke(-1)$, the resulting homomorphism $V^\ast\otimes
H^1\kf(-1)\to H^1\kf$ is described by the matrix
\[
\left(
  \begin{array}{c|c}
v_1 v_2 & 0\\
v_3 v_4 \\\hline
0 & v_1 v_2\\
  & v_3 v_4
  \end{array}
\right)
\]
which describes at the same time the right part $4 \Omega^1(1)\to 4\ko$
of the Beilinson II monad of $\kf$. It follows that $h^1\kf(1)=0$ in this case.

\vskip1cm

\section{Irreducibility of $MI(5)$}\label{section7}

In this section the irreducibility of $MI(5)$ is proved. The proof is mainly
based on the properties of the sets $M(4,4)_\xi$ and the fact that for any
$\ke_\omega\in MI(5)$ there is a $\xi\in H^1\ke_\omega(-1)$ such that
$(4,4)$--instanton $\ke_{\bar{\omega}}$ satisfies
$h^1\ke_{\bar{\omega}}(1)\leq 1$. The {\bf plan of the proof} is the
following. Let us recall the notations
\[
M(4,4)=\{\omega\in S^2H_4^\ast\otimes \Lambda^2 V^\ast\ |\ \rk(\omega)=12,\
\omega \text{ non--degenerate}\}
\]
and
\[
M^0(4,4)=\{\omega\in M(4,4)\ |\ H^1\ke_\omega(1)=0\}.
\]
In addition we need a partial completion of $M(4,4)$, by forgetting the
non--degeneracy,
\[
\widetilde{M}(4,4): = \{\omega\in S^2H_4^\ast\otimes \Lambda^2 V^\ast\ |\
\rk(\omega)=12\}
\]
and for $\xi\in H^\ast_4$, consistent with the earlier notation,
\[
\widetilde{M}(4,4)_\xi:=\{\omega\in\widetilde{M}(4,4)|\quad 
\Im(\omega)\cap(\xi\otimes V^\ast)=0\}.
\]
We have 
\[
\begin{array}{lcl}
M(4,4)_\xi & = & M(4,4)\cap\widetilde{M}(4,4)_\xi\\
M^0(4,4)_\xi & = & M^0(4,4)\cap\widetilde{M}(4,4)_\xi
\end{array}
\]
The space $\widetilde{M}(4,4)$ is introduced for technical reasons which
become apparent in Lemma \ref{affbdl}. For an element
$\omega\in\widetilde{M}(4,4)$ we only get a sheaf $\ke_\omega$ from the monad
  construction \ref{monc}. There is the morphism
\[
\rho=\res_\xi:\widetilde{M}(4,4)_\xi\to M(3,6)=M^0(3,6)
\]
assigning to $\omega$ the map $\bar{H}_3\otimes
V\xrightarrow{\bar{\omega}}\bar{H}_3^\ast\otimes V^\ast$, where
$\bar{H}_3=\Ker(\xi)$. Because $\rk(\bar{\omega})=12$, this map is an
isomorphism and hence non--degenerate. We are going to prove
\vskip5mm

\begin{sub}\label{affbdl} {\bf Lemma}: $\widetilde{M}(4,4)_\xi\xrightarrow{\rho}M(3,6)$ is
  an affine bundle of dimension $54$ and of fibre dimension $18$.
\end{sub}

Because $M(3,6)$ is an open part of $S^2\bar{H}_3^\ast\otimes \Lambda^2
V^\ast$, it is smooth and irreducible of dimension $36$, and therefore
$\widetilde{M}(4,4)_\xi$ is smooth and irreducible of dimension $54$.
Then also the open subset $M(4,4)_\xi$ of $\widetilde{M}(4,4)_\xi$ is smooth
and irreducible of dimension $54$.
\vskip5mm

{\bf Remark}: Irreducibility and smoothness of $M(4,4)_\xi$ is already proved
\ref{rm44}.  The bundle structure and the partial
completion will be used to prove

\begin{sub}\label{codim2} {\bf Lemma}: For any $\xi\in
  H_4^\ast\smallsetminus\{0\},\,  \dim (M(4,4)_\xi\smallsetminus
  M^0(4,4)_\xi)\leq \dim M(4,4)_\xi-2$.
\end{sub}
\vskip5mm

Assuming the two lemmata, the irreducibility of $MI(5)$ is achieved with the
following arguments. The open part
\[
W:=\underset{\eta\neq 0}{\cup} M(4,4)_\eta
\]
of $M(4,4)$ is smooth and irreducible by \ref{union}. We let 
\[
W^0:=\{\omega\in W\ |\ H^1\ke_\omega(1)=0\}=\underset{\eta\neq 0}{\cup}
M^0(4,4)_\eta.
\]
and
\[
W^1:=\{\omega\in W\ |\ h^1\ke_\omega(1)\leq 1\}.
\]
Then $W^0\subset W^1\subset W$ are open subsets. It follows from Lemma
\ref{codim2}, that 
$$
\codim(W^1\smallsetminus W^0)\leq \dim W^1-2.\eqno(1)
$$

Now we consider the morphisms $\omega\mapsto \bar{\omega}=\res_\xi \omega$
\[
M(5,2)_\xi\to M(4,4)
\]
and the inverse images under $r=\res_\xi$,
\[
\begin{array}{ccc}
U_\xi^1 & \xrightarrow{r_1} & W^1\\
\cup & & \cup\\
U_\xi^0 & \xrightarrow{r_0} & W^0.
\end{array}
\]
By Proposition \ref{boundh1}, Lemma \ref{rkl}, and Proposition \ref{exxi}, any element
$\omega\in M(5,2)$ is contained in one of the open sets $U_\xi^1$, i.e.
\[
M(5,2)=\underset{\xi\neq 0}{\cup} U_\xi^1,
\]
but we don't know whether the same holds true for the sets $U_\xi^0$. By
\ref{fibresxi} the fibres of $r_1$ and $r_0$ are open sets of a linear space
and their dimension are
\[
\begin{array}{lcl}
\dim r_1^{-1}(\bar\omega) & = & 4+h^0\ke_{\bar{\omega}}(1)\leq 9\\
\dim r_0^{-1}(\bar\omega) & = & 4+h^0\ke_{\bar{\omega}}(1)=8.
\end{array}
\]
Because the difference of the fibre dimensions is at most $1$, the estimate
(1) implies
$$
\dim(U_\xi^1\smallsetminus U_\xi^0)<\dim U_\xi^0\eqno (2)
$$
for any $\xi\neq 0$. It follows that
\[
\dim(M(5,2)\smallsetminus\underset{\xi\neq 0}{\cup}
U_\xi^0)<\dim\underset{\xi\neq 0}{\cup} U_\xi^0=62.
\]
Because any component of $M(5,2)$ has dimension $\geq 25+40-3=62$ and because
$\underset{\xi\neq 0}{\cup} U_\xi^0$ is irreducible, it follows that $M(5,2)$
is irreducible of the expected dimension $62$.

It follows now from Lemma \ref{crirr} that also $MI(5)$ is irreducible of the
expected dimension $37$. We thus have
\vskip5mm

\begin{sub}\label{thmirr} {\bf Theorem}: $MI(5)$ is irreducible of dimension $37$.
\end{sub}
\vskip5mm

\begin{sub}\label{prf71} \rm {\bf Proof of Lemma \ref{affbdl}:}
We show in fact that 
\[
\widetilde{M}(4,4)_\xi\cong M(3,6)\times\Hom_k(H_3, \Lambda^2 V^\ast).
\]
For the proof we need to distinguish between the linear maps
$A\xrightarrow{\varphi} B\otimes \Lambda^2 V^\ast$ and the corresponding
operators $A \otimes V\xrightarrow{\widetilde{\varphi}} B\otimes V^\ast$
which are skew with respect to $V$ for any two vector spaces $A$ and $B$. 
Because $\xi$ is fixed, we can choose a
decomposition
\[
H_1\oplus \bar{H}_3=H_4,
\]
where $\bar{H}_3$ is the kernel of $\xi$. Given $\bar{\omega}\in M(3,6)$ and
$\bar{H}_3\xrightarrow{\alpha} \Lambda^2 V^\ast\cong 
H_1^\ast\otimes\Lambda^2 V^\ast$, we obtain the operator
\[
\widetilde{\omega}=\left(
  \begin{array}{c|c}
\widetilde{\alpha}\widetilde{\bar{\omega}}^{-1}\widetilde{\alpha}^\ast &
\widetilde{\alpha}\\\hline
\widetilde{\alpha}^\ast & \widetilde{\bar{\omega}}
  \end{array}\right):(H_1\oplus\bar{H}_3)\otimes V\to (H_1^\ast\oplus
\bar{H}^\ast_3)\otimes V^\ast
\]
where $\alpha^\ast$ denotes the dual of $\alpha$ with respect to $\bar{H}_3$ and
$H_1$. Because the upper row of $\widetilde{\omega}$ is a combination of the
lower row, we have
$\rk(\widetilde{\omega})=\rk(\widetilde{\bar{\omega}})=12$. It is clear that
$\widetilde{\omega}$ is skew with respect to $V$. Therefore
$\widetilde{\omega}$ defines an element $\omega\in\widetilde{M}(4,4)_\xi$. We
thus have a morphism
\[
M(3,6)\times \Hom(\bar{H}_3, \Lambda^2 V^\ast)\to\widetilde{M}(4,4)_\xi.
\]
This is even an isomorphism over $M(3,6)$, because if 
\[
\omega=\left(\begin{smallmatrix}\varphi &
    \alpha\\\alpha^\ast&\bar{\omega}\end{smallmatrix}\right):(H_1\oplus
    \bar{H_3})\to (H_1^\ast\oplus\bar{H}_3^\ast)\otimes \Lambda^2 V^\ast
\]
is in $M(4,4)_\xi$, we have
$12=\rk(\widetilde{\omega})=\rk(\widetilde{\bar{\omega}})$ and it follows that 
\[
(\widetilde{\varphi},\widetilde{\alpha})=\widetilde{\alpha}\widetilde{\bar{\omega}}^{-1}(\widetilde{\alpha}^\ast,\widetilde{\bar{\omega}}).
\]

\end{sub}
\vskip5mm

\begin{sub}\label{prf72}\rm {\bf Proof of Lemma \ref{codim2}:} 
For fixed $\eta\in H_4^\ast\smallsetminus\{0\}$ there is an isomorphism
$$ \widetilde{M}(4,4)_\eta\cong M(3,6)\times\Hom(H_3,
\Lambda^2 V^\ast)=: \widetilde{X}$$ 
by Lemma \ref{affbdl}. To each pair
  $(\bar{\omega}, \alpha)\in\widetilde{X}$ we have the simplified Beilinson II
  presentation of $\ke_{\bar{\omega}}$ together with a homomorphism
  $\sigma(\alpha)$ induced by the diagram
\[
\xymatrix{0\ar[r] & H_3\otimes \ko(-1)\ar[ddr]_\varepsilon\ar[r]^{a(\bar{\omega})} &
    H_3^\ast\otimes \Omega^1(1)\ar[r]\ar@{_{(}->}_j[d] &
    \ke_{\bar{\omega}}\ar[r]\ar@{-->}[dd]^{\sigma(\alpha)} & 0\\ &&
    H_3^\ast\otimes V^\ast\otimes \ko \\
& & H_3\otimes V\otimes
    \ko\ar[u]^\approx_{\bar{\omega}}\ar[r]^-{\widetilde{\alpha}} & \ko(1)\ .
}
\]
Note here that $\widetilde{\alpha}\circ\varepsilon=0$ because $\alpha$ is skew
with respect to $V$, such that $\widetilde{\alpha}\circ\bar{\omega}^{-1}\circ
j$ factors through $\ke_{\bar{\omega}}$. Note further that in case
$(\bar{\omega},\alpha)$ corresponds to an $\omega\in M(4,4)_\eta$, then
$\ke_\omega$ is the cohomology of the monad
\[
\ko(-1)\xrightarrow{\sigma(\alpha)^\ast}\ke^\ast_{\bar{\omega}}
\xrightarrow[j_{\bar{\omega}}]{\approx}\ke_{\bar{\omega}}\xrightarrow{\sigma(\alpha)}\ko(1)
\]
in which case $\sigma(\alpha)$ is surjective, see beginning of Section \ref{section3}.

We let $X\subset \widetilde{X}$ be the open subset corresponding to
$M(4,4)_\eta$, or defined by the surjectivity of $\sigma(\alpha)$. In addition
we let $\widetilde{X}^0\subset\widetilde{X}$ denote the open part of
$\widetilde{X}$ where 
\[
H^0\sigma(\alpha)(1): H^0\ke_{\bar{\omega}}(1)\to H^0\ko(2)
\]
is surjective. Then 
\[
M^0(4,4)_\eta\cong X^0=X\cap \widetilde{X}^0
\]
under the above isomorphism because the cokernel of $H^0\sigma(\alpha)(1)$ is
then isomorphic to $H^1\ke_\omega(1)$. To prove Lemma \ref{codim2}, it is
sufficient to prove that 
$$
\codim(\widetilde{X}\smallsetminus\widetilde{X}^0, \widetilde{X})\geq
2.\eqno(\ast)
$$
Now $(\ast)$ will follow from the following two statements.

{\bf Claim F}: There exists an $\bar{\omega}\in M(3,6)$ such that for the
fibres,  $\widetilde{X}^0_{\bar{\omega}}\subset \widetilde{X}_{\bar{\omega}}$  of $\widetilde{X}^0\subset\widetilde{X}$ over $\bar{\omega}$,
\[
\codim(\widetilde{X}_{\bar{\omega}}\smallsetminus
\widetilde{X}^0_{\bar{\omega}}, \widetilde{X}_{\bar{\omega}})\geq 2.
\]

{\bf Claim B:} Let $\Sigma\subset M(3,6)$ be the closed subvariety of points
$\bar{\omega}$ for which $\widetilde{X}^0_{\bar{\omega}}=\emptyset$, i.e.\ the
set of points $\bar{\omega}$ for which $H^0\sigma(\alpha)(1)$ is not
surjective for any $\alpha$. Then
\[
\codim(\Sigma, M(3,6))\geq 2.
\]

\vskip5mm

\begin{subsub}\label{prfb} \rm {\bf Proof of claim B:} 

a) In order to incorporate a
  projective curve which doesn't meet $\Sigma$, we enlarge $M(3,6)$ and
  $\Sigma$ as follows. Let
\[
\widetilde{M}(3,6):= M(3,6)\cup M^0(3,4).
\]
Then $\widetilde{M}(3,6)\subset S^2H_3^\ast\otimes \Lambda^2 V^\ast$ consists
of non--degenerate $\bar{\omega}$ of rank $10$ if $\bar{\omega}\in M^0(3,4)$
and of rank $12$ if $\bar{\omega}\in M(3,6)$, and such that
$H^1\ke_{\bar{\omega}}(1)=0$. If $\bar{\omega} \in M^0(3,4)$, the Beilinson II
monad of $\ke_{\bar{\omega}}$ is of the type
\[
0\to H_3\otimes \ko(-1)\xrightarrow{a(\bar{\omega})} H_3^\ast\otimes
  \Omega^1(1)\to Q\otimes \ko\to 0
\]
where $Q$ is the cokernel of $\bar{\omega}$.

If $\kf_{\bar{\omega}}$ is the cokernel of $a(\bar{\omega})$, we have the
exact sequence
\[
0\to \ke_{\bar{\omega}}\to\kf_{\bar{\omega}} \to Q\otimes \ko\to 0.
\]
Given a second component $\alpha\in \Hom(H_3, \Lambda^2 V^\ast)$, we also
obtain a homomorphism $\kf_{\bar{\omega}}\xrightarrow{\sigma(\alpha)} \ko(1)$
by the diagram
\[
\xymatrix{& & H_3\otimes V\otimes
  \ko\ar[d]^{\bar{\omega}}\ar[r]^>>>>{\widetilde{\alpha}} & \ko(1) \\
& &H_3^\ast\otimes V x_\otimes \ko\ar@{-->}[ur] &\\
0\ar[r] & H_3\otimes \ko(-1)\ar[r] & H_3^\ast\otimes
  \Omega^1(1)\ar@{_{(}->}[u]\ar[r] &
  \kf_{\bar{\omega}}\ar@{-->}[uu]_{\sigma(\alpha)}\ar[r] & 0
}
\]
which in case $\bar{\omega}\in M(3,6)$ coincides with
$\ke_{\bar{\omega}}\xrightarrow{\sigma(\alpha)}\ko(1)$. The surjectivity of
$H^0\sigma(\alpha)(1)$ does not depend on the choice of the factorization of
$\widetilde{\alpha}$.

We let now 
$\widetilde{\Sigma}\subset\widetilde{M}(3,6)$
be the locus of points $\bar{\omega}$ for which $H^0\sigma(\alpha)(1)$ is not
surjective for any $\alpha$. By this definition we have
\[
\Sigma=\widetilde{\Sigma}\cap M(3,6).
\]

b) Because $M_{10}(H_3)$ had been shown to be an irreducible hypersurface in
$S^2H_3^\ast\otimes \Lambda^2 V^\ast$  whose complement is $M(3,6)$, see
\ref{m10}, and because $M^0(3,4)$ is an open part of it, the complement of
$\widetilde{M}(3,6)$ in $S^2H_3^\ast\otimes \Lambda^2 V^\ast$ has codimension
$\geq 2$. In order to show that $\codim(\widetilde{\Sigma},
\widetilde{M}(3,6))\geq 2$, we first construct an embedding $k^2\smallsetminus
\{0\}\to\widetilde{M}(3,6)\smallsetminus \widetilde{\Sigma}$. 

c) This embedding
is defined as follows. We let $e_1, \ldots, e_4$ be the standard basis of
$H_4=k^4$ and $e_1^\ast, \ldots, e_4^\ast$ be its dual basis. We let 
\[
H_4\xrightarrow{\omega} H_4^\ast\otimes \Lambda^2 V^\ast
\]
be given by the matrix
\[
\omega=\left(
  \begin{array}{cccc}
\omega'_{11} & \omega'_{12} & 0 & 0\\
\omega'_{12} & \omega'_{22} & 0 & 0\\
0            & 0            & \omega''_{11} & \omega''_{12}\\
0            & 0            & \omega''_{12} & \omega''_{22}
  \end{array}\right)
\]
with $\omega'_{ij},\, \omega''_{ij}\in\Lambda^2 V^\ast$ which represents the
direct sum of two $2$--instantons $\ke'$ and $\ke'',
\ke_\omega=\ke'\oplus\ke''$. Then $H^1\ke_\omega(1)=0$. For $t=(t_0, t_1)\neq
0$ we consider
\[
\xi_t=(-t_1, 0, t_0, 0)\in H_4^\ast
\]
and its kernel
\[
0\to k^3\xrightarrow{f_t} k^4\xrightarrow{\xi_t} k\to 0
\]
defined by the matrix
\[f_t=\left(
  \begin{array}{ccc}
t_0 & 0 & 0\\
0 & 1 & 0\\
t_1 & 0 & 0\\
0 & 0 & 1
  \end{array}\right)\ .
\]
We let $\bar{\omega}_t=(f_t^\ast\otimes\id)\circ \omega\circ f_t$ such that
\[
\bar{\omega}_t = \left(
  \begin{array}{ccc}
t_0^2\omega'_{11}+t_1^2\omega''_{11} & t_0\omega'_{12} & t_1\omega''_{12}\\
t_0\omega'_{12} & \omega'_{22} & 0\\
t_1\omega''_{12} & 0 & \omega''_{22}
  \end{array}\right)\ .
\]
We may assume that both $\omega'_{22}$ and $\omega''_{22}$ have rank $4$ as
operators $V\to V^\ast$. Then by an elementary matrix operation we can kill
$t_0\omega'_{12}$ and $t_1\omega''_{12}$ in the first row and obtain a matrix
\[
\left(
  \begin{array}{ccc}
t_0^2\eta'_{11}+t_1^2\eta''_{11} & 0 & 0\\
t_0\omega'_{12} & \omega'_{22} &\\
t_1\omega''_{12} & & \omega''_{22}
  \end{array}\right)
\]
with $\eta'_{11}=\omega'_{11}-\omega'_{12}\omega_{22}^{\prime-1}\omega'_{12}$ and
$\eta''_{11}=\omega''_{11}- \omega''_{12}\omega_{22}^{''-1}\omega''_{12}$.

Because $\omega'$ and $\omega''$ represent $2$--instantons,
$\rk(\omega')=\rk(\omega'')=6$ and therefore
\[
\rk\ \eta'_{11}=\rk\  \eta''_{11}=2.
\]
But we can choose $\omega'$ and $\omega''$ such that
$\rk\ \omega'_{22}=\rk\ \omega''_{22}=4$ and in addition
\[
\Im\, \eta'_{11}+\Im\ \eta''_{11}=V^\ast.
\]
Then $t_0^2\eta'_{11}+t_1^2\eta''_{11}$ is an isomorphism $V\to V^\ast$ for
$t_0t_1\neq 0$. Therefore, with this choice of $\omega'$ and $\omega''$, we
have
\[
\begin{array}{lcl}
\bar{\omega}_t\in M(3,6)& \text{ for } & t_0t_1\neq 0\\
\bar{\omega}_t\in M^0(3,4) & \text{ if } &  t_0=0\text{ or } t_1=0.
\end{array}
\]
The first statement follows directly from the fact that $\rk\ \bar{\omega}_t=12$
if $t_0t_1\neq 0$. If $t_0=0$ or $t_1=0$, then $\rk\ \bar{\omega}_t=10$. In that
case, e.g.\ $t_1=0$,
\[
\bar{\omega}_t=\left(
  \begin{array}{ccc}
t_0^2\omega'_{11} & t_0\omega'_{12} & 0\\
t_0\omega'_{12} & \omega'_{22} & 0\\
0 & 0 & \omega''_{22}
  \end{array}
\right)
\]
and then the bundle of $\bar{\omega}_t$ is the direct sum
\[
\ke_{\bar{\omega}_t}\cong \ke_{\omega'}\oplus\ke_{\omega''_{12}}
\]
with $\ke'=\ke_{\omega'}$ a $2$--instanton and $\ke_{\omega''_{22}}$ a
$1$--instanton (null--correlation bundle). So $\bar{\omega}_t$ is
non--degenerate and $H^1\ke_{\bar{\omega}_t}(1)=0$. It follows that $t\mapsto
\bar{\omega}_t$ is a morphism
\[
k^2\smallsetminus\{0\}\to\widetilde{M} (3,6)
\]
which is an embedding by the shape of the matrix $\bar{\omega}_t$. Moreover,
$\bar{\omega}_t\not\in\widetilde{\Sigma}$ for any $t$. To see this, we
consider the case $t_1\neq 0$ first. In that case the linear embedding
$k^3\xrightarrow{f_t} k^4$ can be considered as the kernel of a fixed
$\xi=(1,0,0,0)$ by the following diagram
\[
\xymatrix{0\ar[r] & k^3\ar@{=}[d]\ar[r]^{f_t} &
  k^4\ar[d]_\approx^{a_t}\ar[r]^{\xi_t} & k\ar@{=}[d]\ar[r] & 0\\
0\ar[r] & k^3\ar[r]^g & k^4\ar[r]^\xi & k\ar[r] & 0
}
\]
where $\xi_t=(-t_1, 0, t_0, 0)$ and $\xi=(1,0,0,0)$ with
\[
g=\left(
  \begin{array}{ccc}
0 & 0 & 0\\
1 &   &\\
& 1 &\\
& & 1
  \end{array}\right)\quad\text{ and }\quad a_t=\left(
  \begin{array}{cccc}
-1 & 0 & t_0/t_1 & 0\\
0 & 0 & 1 & 0\\
0 & 1 & 0 & 0\\
0 & 0 & 0 & 1
  \end{array}\right)\raisebox{-4ex}{.}
\]
If we denote 
\[
\omega_t:=a_t^{-1,\ast}\omega a_t^{-1}
\]
we have
\[
\bar{\omega}_t = g^\ast \omega_t g=\res_\xi\omega_t.
\]
Because $\omega_t\sim\omega$, we have $\ke_{\omega_t}\cong\ke_\omega$ and so
$h^1\ke_{\omega_t}(1)=0$. If $\bar{\omega}_t\in M(3,6)$, then
$\bar{\omega}_t\not\in \Sigma$. If, however, $\bar{\omega}_t\in M^0(3,4)$,
then $\ke_{\bar{\omega}_t}$ is the cohomology of the complex
\[
\xymatrix{\ko(-1)\ar@/_1pc/[rr] \ar[r]& \kf_{\bar{\omega}_t}^\ast \ar[r] &
  \kf_{\bar{\omega}_t}\ar[r]& \ko(1)}.
\]
\vskip3mm

Then also in this case $H^0\kf_{\bar{\omega}_t}(1)\to H^0\ko(2)$ has cokernel
$H^1\ke_{\omega_t}(1)=0$. Therefore $\bar{\omega}_t\not\in\widetilde{\Sigma}$
for $t_1\neq 0$. By symmetry in $t_0, t_1$ , also
$\bar{\omega}_t\not\in\widetilde{\Sigma}$ for $t_0\neq 0$. Therefore we have
an embedding
\[
k^2\smallsetminus \{0\}\hookrightarrow\widetilde{M}(3,6)\smallsetminus
\widetilde{\Sigma}.
\]

d) The map $t\mapsto \bar{\omega}_t$ cannot directly be used to define a line
in $\P S^2H_3^\ast\otimes \Lambda^2 V^\ast$ because $\bar{\omega}_t$ is not
homogeneous in $t$. But it can be used to construct a projective line in
$\P H_4^\ast$ which doesn't intersect a transformation of $\widetilde{\Sigma}$
so that we can conclude from that that $\widetilde{\Sigma}$ has codimension
$\geq 2$. To do this we consider the following transformation diagram
\[
\P H_4^\ast\xleftarrow{\pi} \P\Hom(H_3, H_4)^0\xrightarrow{\theta}
\P S^2H_3^\ast \otimes\Lambda^2 V^\ast.
\]
In this diagram $\P\Hom(H_3, H_4)^0$ denotes that open set of injective maps
and $\pi$ is the principal $PGL(H_3)$--bundle over the Grassmannian
$G(3,H_4)=\P H_4^\ast$. The morphism $\theta$ is defined by
\[
f\overset{\theta}{\mapsto}(f^\ast\otimes \id)\circ \omega\circ(f\otimes \id)
\]
using the form of $\omega$ in c). We let
\[
\widehat{\Sigma}\subset\widehat{M}(3,6)\subset\P S^2 H_3\ast\otimes \Lambda^2
V^\ast
\]
denote the subvarieties of which $\widetilde{\Sigma}$ and $\widetilde{M}(3,4)$
are the affine cones. So $\widehat{\Sigma}$ is closed in the open set
$\widehat{M}(3,6)$ of codimension $\geq 2$. Now
\[
\theta^{-1}\widehat{\Sigma}\subset\theta^{-1}\widehat{M}(3,6)\subset\P\Hom(H_3,
H_4)^0
\]
are $PGL(H_3)$--invariant and therefore there are an open subset $U\subset \P
H_4^\ast$ and a closed subscheme $Z\subset U$ such that 
\[
\theta^{-1}\widehat{\Sigma}=\pi^{-1} Z\quad \text{ and }\quad
\theta^{-1}\widehat{M}(3,6)=\pi^{-1} U.
\]
Moreover, there is the irreducible hypersurface
$\widehat{M}_{10}(H_3)\subset\P S^2H_3^\ast\otimes \Lambda^2V^\ast$ defined by
$\rk\ \bar{\omega}\leq 10$. Its inverse image $\theta^{-1}\widehat{M}_{10}
(H_3)$ is the subvariety of those $\langle f\rangle \in\P\Hom(H_3, H_4)^0$,
for which $(H_4/\Im(f))^\ast\otimes V^\ast\to H^1\ke_\omega$ is not an
isomorphism, see \ref{rkl}, because then $\bar{\omega}=\theta(f)$ has rank
$10$. It is then easily seen that this condition defines an effective divisor,
which is $\theta^{-1}\widehat{M}_{10}(H_3)$. Because this is also
$PGL(H_3)$--invariant, there is a divisor $D\subset \P H_4^\ast$ such that
\[
\theta^{-1}\widehat{M}_{10}(H_3)=\pi^{-1} D.
\]
Now the family $(f_t)$ of c) defines a linear embedding $\langle t\rangle \mapsto
\pi(\langle f_t\rangle )=\langle \xi_t\rangle $
\[
\P_1\xrightarrow[\approx]{} L\subset \P H_4^\ast
\]
and by the result of c) we have 
\[
L\subset U\smallsetminus Z\quad \text{ and }\quad L\not\subset D.
\]
Let $Y\subset \P H_4^\ast$ be the complement of $U$. Then $L$ doesn't
intersect $Y\cup Z=Y\cup \bar{Z}$ and so $Z$ is of codimension $\geq 2$. Then
also $\theta^{-1}\widehat{\Sigma}$ has codimension $\geq 2$ in
$\theta^{-1}\widehat{M}(3,6)$. But since also the complement of
$\theta^{-1}\widehat{M}(3,6)$ has codimension $\geq 2$ in $\P\Hom(H_3, H_4)$,
there is a line $L'\subset \P\Hom(H_3, H_4)$ such that
\[
L'\subset\theta^{-1}\widehat{M}(3,6)\smallsetminus
\theta^{-1}\widehat{M}\quad\text{ and }\quad
L'\not\subset\theta^{-1}\widehat{M}_{10}(H_3).
\]
We let  now $\Gamma=\theta(L')$. This is a complete curve in
$\P S^2H_3^\ast\otimes S^2 V^\ast$ (not contracted to a point because $L'$
intersects the divisor $\theta^{-1}\widehat{M}_{10}(H_3)$ and is not contained
in it) and such that
\[
\Gamma\subset \widehat{M}(3,6)\smallsetminus \widehat{\Sigma}.
\]
Then $\Gamma$ cannot intersect the closure of $\widehat{\Sigma}$ in $\P
S^2H_3^\ast\otimes \Lambda^2 V^\ast$. So we have shown that $\codim
(\widehat{\Sigma}, \widehat{M}(3,6))\geq 2$. This is equivalent to claim B.
\end{subsub}

\begin{subsub}\label{rmnc} \rm {\bf Remarks on null--correlation bundles}:

The tensors $\bar{\omega}$ provided for claim F are diagonal matrices of
elements $ \eta\in\Lambda^2 V^\ast$ of rank $4$ (or indecomposable) which
correspond to $1$--instantons or null-correlation bundles (nc--bundles for
short). Given $\eta$ of rank $4$, we have the defining sequence $0\to
\ko(-1)\to\Omega^1(1)\to\ke_\eta\to 0$. Then for any $\alpha\in \Lambda^2 V^\ast$
we obtain a homomorphism $\ke_\eta\xrightarrow{\sigma(\alpha)}\ko(1)$ as in
\ref{prfb}, using the diagram
\[
\xymatrix{0\ar[r] & \ko(-1)\ar[r]^\eta & \Omega^1(1)\ar[r]\ar@{_{(}->}[d] &
  \ke_\eta\ar[r]\ar@{-->}[dd]^{\sigma(\alpha)} & 0\\
 & & V^\ast\otimes\ko & &\\
 & & V\otimes \ko\ar[u]^\approx_\eta\ar[r]^{\widetilde{\alpha}} & \ko(1) & 
}.
\]
One can easily see that $\sigma(\alpha)=0$ if and only if $\alpha$ and $\eta$
are proportional in $\Lambda^2 V^\ast$ or $\langle \alpha\rangle =
\langle \eta\rangle $. Therefore, a
non--zero image $\ki$ of $\ke_\eta(-1)\to\ko$ is determined by a line
$g=\overline{\langle \alpha\rangle , \langle \eta\rangle }$ in $\P\Lambda^2 V^\ast$ and is independent of
$\langle \alpha\rangle $ on $g$. We therefore write $\ki_g$ for the image of $\ke_\eta(-1)$
under $\sigma(\alpha)$.

In this section we identify the Grassmannians $G=G(2,V)$ and $G(2,V^\ast)$ via
the Pl\"ucker embeddings in $\P\Lambda^2 V\cong P\Lambda^2 V^\ast$, because
$\Lambda^2 V\cong\Lambda^2 V^\ast\otimes \Lambda^4 V$ canonically. Then for a
line $l$ in $\P V$ through $\langle x\rangle , \langle y\rangle $ with equations $z, w$ we write
$l=\langle x\wedge y\rangle =\langle z\wedge w\rangle $.
\vskip5mm

 {\bf Lemma 1}: {\em Let $\langle \eta\rangle \in\P\Lambda^2
  V^\ast\smallsetminus G$ and let $g$ be a line in $\P \Lambda^2 V^\ast$
  through $\langle \eta\rangle $ with associated ideal sheaf $\ki_g$. Then
  \begin{enumerate}
  \item [a)] If $g$ meets the Grassmannian $G$ in two points $l_1\neq l_2$,
    then $l_1\cap l_2=\emptyset$ in $\P V$ and $\ki_g= \ki_{l_1\cup l_2}$.
\item [b)] If $g$ is tangent to $G$ at $l$, then $\ki_g$ is the twisted double
  structure on the line $l$ given by the tangent direction along $g$.
 \end{enumerate}
In both cases $h^0\ki_g(2)=4.$}

\begin{proof} Let $Z$ be the zero scheme of $\ki_g$. We have the exact
  sequence $0\to \ko(-1)\to\ke_\eta\to\ko(1)\to\ko_Z(1)\to 0$ and by that for
  the Hilbert polynomial $\chi\ko_Z(m)=2m+2$. This means that $Z$ is a pair of
  disjoint lines or a twisted double structure on a line. In order to prove
  that $\ki_g$ is the ideal of the lemma, it is sufficient to determine the
  space $H^0\ki_g(2)\subset S^2 V^\ast$ by the induced diagram
\[
\xymatrix{0\ar[r] & k\ar[r]^\eta & \Lambda^2 V^\ast\ar[r]\ar@{_{(}->}[d] &
  H^0\ke_\eta(1)\ar[r]\ar[dd]^{H^0\sigma(\alpha)(1)} & 0\\
& & V^\ast\otimes V^\ast & &\\
& & V\otimes V^\ast\ar[u]^\approx_{\eta\otimes\id}\ar[r]^{\widetilde{\alpha}} &
  S^2 V^\ast & }
\]
where $g\ni\langle \alpha\rangle \neq \langle \eta\rangle $.

{\bf case a)} Because $\langle \alpha\rangle , \langle \eta\rangle , 
l_1, l_2$ are on the same line
$g\not\subset G$, we may choose a basis $e_0, \ldots, e_3\in V$ with dual basis
$z_0, \ldots z_3$ such that $l_1=\langle e_0\wedge e_1\rangle$, 
$l_2=\langle e_2\wedge e_3\rangle$, $\eta=z_0\wedge z_1+z_2\wedge z_3$ and 
$\alpha=\lambda z_0\wedge z_1+\mu z_2\wedge
z_3$ with $\lambda\neq \mu$. With this choice one computes that in the above
diagram the image of $H^0\sigma(\alpha)(1)$ is spanned by the 4 forms $z_0z_2,
z_0z_3, z_1 z_2, z_1 z_3$ which are the generators of the ideal of $l_1\cup
l_2$.

{\bf case b)} Here we may choose the same $\eta$ but $\alpha=z_1\wedge z_3$
such that $\langle \alpha\rangle =\langle e_0\wedge e_2\rangle =l$ in $\P\Lambda^2 V^\ast=\P\Lambda^2
V$. With the same type of calculation we find that $H^0\ki_g(2)$ is spanned by
the 4 forms $z_1^2, z_1z_3, z_3^2, z_0z_1+z_2z_3$, which are the generators of
the ideal of the double structure on $l$.
\end{proof}

{\bf Lemma 2:} {\em (1) Let $p_1, p_2$ be two different points on a line $l$
  in $\P_3$ and let $l_1, l_2$ be two skew lines not meeting $p_1, p_2$. Then
  the space of quadrics containing $p_1, p_2, l_2, l_2$ is $2$--dimensional if
  and only if $l$ doesn't meet both $l_1$ and $l_2$, and otherwise this space
  is $3$--dimensional.

(2) If $l_0$ is a twisted double line as in Lemma 1 above, which doesn't meet
$l$, then the space of quadrics containing $p_1, p_2$ and $l_0$ is
$2$--dimensional, too.}

\begin{proof} By elementary computation using suitable coordinates.
\end{proof}
\end{subsub}
\vskip5mm
\begin{subsub}\label{3ncb} {\bf Three nc--bundles}:\rm

We consider first two nc--bundles $\ke_1=\ke_{\eta_1}, \ke_2=\ke_{\eta_2}$
with homomorphisms $\ke_i\xrightarrow{\sigma_i} \ko(1)$, and ideal sheaves
$\ki_1, \ki_2$ respectively. Then the image of $H^0(\ke_1(1)\oplus\ke_2(1))\to
H^0\ko(2)$ under $\sigma_1+\sigma_2$ is $H^0\ki_1(2)+H^0\ki_2(2)$, and
similarly for three such data $(\ki_i, \sigma_i)$. We let $Z_1, Z_2$ be the
zero schemes of $\ki_1, \ki_2$ respectively.

{\bf Lemma 3}: {\em Let $Z_1=l_{11}\cup l_{12}$ and $Z_2=l_{21}\cup l_{22}$
  with all 4 lines different. Suppose that either\\
(i)  $Z_1\cap Z_2=\emptyset$ and the 4 lines do not belong to the same
  ruling of a quadric, or\\
(ii) $Z_1\cap Z_2=\{x_1, x_2\}\subset l_{11}$ and 
$Z_2\cap l_{12}=\emptyset$.\\
Then
\[
H^0\ki_1(2)\cap H^0\ki_2(2)=0.
\]}

\begin{proof} (i) well--known. (ii) Let $Q$ be a smooth quadric containing
  $l_{11}, l_{12}$ (in one ruling). Because $l_{2\nu}\cap l_{12}=\emptyset$ and
  $l_{2\nu}\cap l_{11}=\{x_\nu\}$, we have $l_{21}, l_{22}\not\subset
  Q$. Therefore the intersection of the $H^0\ki_\nu(2)$ does not contain a
  smooth quadric. Then any $Q\in H^0\ki_1(2)\cap H^0\ki_2(2)$ is not smooth and 
contains the skew lines $l_{11}$ and $l_{12}$. So 
$Q=P_1\cup P_2$ is a pair of planes. Let $l_{11}\subset
  P_1$. Then $l_{12}\not\subset P_1$ because $l_{11}\cap l_{12}=\emptyset$ and so
  $l_{12}\subset P_2$. Then $l_{21}, l_{22}\not\subset P_2$ and $l_{21},
  l_{22}\subset P_1$. But $l_{21}\cap l_{22}=\emptyset$, a contradiction.
\end{proof}

{\bf Remark}: In case (ii) the $8$--dimensional space
$H^0\ki_1(2)+H^0\ki_2(2)$ equals $H^0\ki_{\{x_1, x_2\}}(2)$.
\vskip5mm

{\bf Corollary}: {\em Let the situation be as in Lemma 3 and let
$\ke_3\xrightarrow{\sigma_3}\ko(1)$ be a third homomorphism of an nc--bundle
with ideal $\ki_3$ and zero scheme $Z_3$. If $Z_3$ consists of two (then skew)
lines $l_{31}, l_{32}$ or of a twisted double line $l_3$, such that the line
$$l_{11}=\overline{x_1,x_2}$$ doesn't meet both of $l_{31}$and $l_{32}$, then
\[
H^0\ki_1(2)+H^0\ki_2(2)+H^0\ki_3(2)= H^0\ko(2).
\]}

\begin{proof} By Lemma 3 the sum of the first two is $H^0\ki_{\{x_1, x_2\}}
  (2)$. By Lemma 2 of \ref{rmnc},
\[
\dim H^0\ki_{\{x_1, x_2\}}(2)\cap H^0\ki_3(2)=2.
\]
\end{proof}
\end{subsub}

\begin{subsub}\label{prff} {\bf Proof of claim F:}\rm 

a) We choose an isomorphism $\bar{H}_3\cong k^3$ and
\[
\bar{\omega}=\left(
  \begin{array}{ccc}
\eta_1 & 0 & 0\\
0 & \eta_2 & 0\\
0 & 0 & \eta_3
  \end{array}\right)
\]
where $\eta_i\in\Lambda^2 V^\ast$ are indecomposable. Let $\ke_i=\ke_{\eta_i}$
be the associated nc--bundles such that
\[
\ke=\ke_{\bar{\omega}}=\ke_1\oplus\ke_2\oplus \ke_3.
\]
By the definition of $\widetilde{X}_{\bar{\omega}}$ in \ref{prf72} this fibre
is isomorphic to $\Hom(k^3, \Lambda^2 V^\ast)$ and we have a surjective linear
map
\[
\widetilde{X}_{\bar{\omega}}\to\Hom(\ke_1\oplus\ke_2\oplus\ke_3, \ko(1))
\]
by
\[
(\alpha_1, \alpha_2, \alpha_3)\mapsto (\sigma_1(\alpha_1),\sigma_2(\alpha_2),
\sigma_3(\alpha_3)).
\]
The open set $\widetilde{X}^0_{\bar{\omega}}$ in
$\widetilde{X}_{\bar{\omega}}$ is then the inverse image of the set of those 
$(\sigma_1, \sigma_2, \sigma_3),\ \ke_i\xrightarrow{\sigma_i} \ko(1)$, such that
for the corresponding ideal sheaves $\ki_i$ we have
$$
H^0\ki_1(2)+H^0\ki_2(2)+H^0\ki_3(2)=H^0\ko(2).\eqno(\ast)
$$
We put $W_i=\Hom(\ke_i,\ko(1))$ and let 
\[
Z\subset\P(W_1\oplus W_2\oplus W_2)
\]
be the closed subvariety of points $\langle \sigma_1, \sigma_2, \sigma_3\rangle $ for which
$(\ast)$ is not satisfied. Then $\codim\ Z\geq 2$ implies claim F for the
chosen $\bar{\omega}$.

b) In order to prove $\codim Z\geq 2$ we consider the natural projection
\[
\P(W_1 \oplus W_2\oplus W_3)\smallsetminus Z'\xrightarrow{\pi}\P W_1\times
\P W_2\times \P W_3
\]
where $Z'$ is the subvariety of points $\langle \sigma_1, \sigma_2, \sigma_3\rangle $ with
at least one component equal to $0$. It is easy to see that $\pi$ is a
principal bundle with fibre and group $(k^\ast)^3/k^\ast$. We have $Z'\subset Z$ and
$Z\smallsetminus Z'$ is invariant under $(k^\ast)^3/k^\ast$. Then
$Y=\pi(Z\smallsetminus Z')$ is closed and $Z\smallsetminus
Z'=\pi^{-1}(Y)$. Now $Y$ is the subvariety of triples $(\langle \sigma_1\rangle , 
\langle \sigma_2\rangle,
\langle \sigma_3\rangle )$ for which $(\ast)$ is not satisfied for the images
$\ki_\nu(1)=\Im\ \sigma_\nu$. Claim F will be proved if $\codim\ Y\geq 2$. This
follows now from

c) {\bf Proposition}: {\em Let $\langle \eta_1\rangle , \langle \eta_2\rangle , 
\langle \eta_3\rangle \in \P\Lambda^2
  V^\ast\smallsetminus G$ be in general position (not co-linear). Then the
  subvariety $Y$ has codimension $\geq 2$.}

\begin{proof} We let $y_i$ denote the points $\langle \eta_i\rangle $ in $\P\Lambda^2 V$.

(i) Let $H(y_2)$ be the polar hyperplane of $y_2$ with respect to the quadric
$G$. There is point $l\in H(y_2)\cap G$ such that the tangent hyperplane $T_l
G$ doesn't contain $y_3$. Denote by 
\[
C(l)=G\cap T_l(G)
\]
the cone of lines in $\P V$ meeting $l$. Then $C(l)\cap H(y_2)$ has
codimension $2$.

(ii) Choose any $l_{11}\in H(y_2)\cap G$ with $y_3\not\in T_{l_{11}} G$ and
such that the line $g_1=\overline{y_1, l_{11}}$ meets $G$ in two different
points $l_{12}\neq l_{11}$. We have $l_{11}, y_2\in T_{l_{11}} G$ because
$l_{11}\in H(y_2)$. Next we choose a line $g_2$ in $T_{l_{11}} G$ through
$y_2$ which meets $G$ in two different points $l_{21}, l_{22}$. These lines
belong to the cone $C(l_{11})$ and thus each meets $l_{11}$ in a point $x_1$
and $x_2$. We may assume that $g_2\cap C(l_{12})=\emptyset$, such that $l_{21}$ and
$l_{22}$ don't meet $l_{12}$. By this choice the conditions of \ref{3ncb}, Lemma
3, are satisfied, and thus 
\[
H^0\ki_{g_1}(2) + H^0\ki_{g_2}(2)= H^0\ki_{\{x_1, x_2\}}(2)
\]
is $8$--dimensional, where $\ki_{g_\nu}$ is the ideal corresponding to the
  line $g_\nu$ through $y_\nu$.

(iii) Let now $P(x_\nu)\subset G$ be the $\alpha$--plane of all lines through
$x_\nu$ and let 
\[
S=P(x_1)\cup P(x_2)\cup (C(l_{11})\cap H(y_3)).
\]
Because $y_3\not\in T_{l_{11}} G,\ C(l_{11})\cap H(y_3)$ is $2$--dimensional
and hence $S$ is $2$--dimensional. Let
\[
\P\Lambda^2 V\smallsetminus \{y_3\}\xrightarrow{\pi} \P(\Lambda^2 V/y_3)
\]
be the central projection. We have $y_3\not\in S$. Then $\pi(S)$ has dimension
$2$ and there is a line $L\subset \P(\Lambda^2 V/y_3)$ which doesn't meet
$\pi(S)$. Let $a_3\in L$ be any point and let $g_3$ denote the line through
$y_3$ in $\P\Lambda^2 V$ given by $\pi^{-1}(a_3)$. If $a_3\not\in
\pi(H(y_3)\cap G)$, the branch locus of $\pi|G$, then $g_3\cap G$ consists of
two different points $l_{31}, l_{32}$ which don't meet $x_1, x_2$. Suppose
that $l_{31}$ meets $l_{11}$, i.e.\ $l_{31}\in C(l_{11})\subset T_{l_{11}} G$.
Then $l_{32}\not\in T_{l_{11}} G$ because otherwise also $y_3\in T_{l_{11}}G$,
which had been excluded by the choice of $l_{11}$. By the corollary in
\ref{3ncb} $H^0\ki_{g_3}(2)$ intersects $H^0\ki_{\{x_1, x_2\}}(2)$ in
dimension $2$. Then for $\langle \alpha_\nu\rangle \in g_\nu, \nu=1, 2, 3$, we
have a point $(\langle \alpha_1\rangle , \langle \alpha_2\rangle , \langle
\alpha_3\rangle )\in\P W_1\times \P W_2\times \P W_3\smallsetminus Y$. If
however, $a_3\in \pi(H(y_3)\cap G)$, then $g_3$ is tangent to $G$ at some
$l_0$ and $\ki_{g_3}$ is a twisted double structure on $l_0$. In this case
$l_0\not\in C(l_{11})$ because $l_0\in H(y_3)$, and so $l_0\cap
l_{11}=\emptyset$. Again by the corollary $H^0\ki_{g_3}(2)$ and
$H^0\ki_{\{x_1, x_2\}}(2)$ intersect in dimension $2$.

(iv) Let now $a'_3, a''_3\in L$ be any two points, $g'_3, g''_3$ the lines in
$\P \Lambda^2 V$ through $y_3$, defined by $\pi^{-1}(a'_3), \pi^{-1}(a''_3)$ and let $\langle \alpha_1\rangle \in g_1 \smallsetminus
\{y_1\}, \langle \alpha_2\rangle \in g_2\smallsetminus \{y_2\}, \langle \alpha'_3\rangle \in
g'_3\smallsetminus \{y_3\}, \langle \alpha''_3\rangle \in g''_3\smallsetminus\{y_3\}$. Then
the points
\[
(\langle \alpha_1\rangle , \langle \alpha_2\rangle , 
\langle s\alpha'_3+t\alpha''_3\rangle )\in\P W_1\times \P W_2\times
\P W_3
\]
do not belong to $Y$ for any $(s,t)\neq (0,0)$. We thus have a complete curve
in the product space not meeting $Y$, i.e.\ $\codim Y\geq 2$. This completes
the proof of claim F and with that the proof of Lemma \ref{codim2} and
finally the proof of theorem \ref{thmirr}.
\end{proof}
\end{subsub}
\end{sub}
\vskip1cm

\section{Smoothness of $MI(5)$}\label{section8}

The induction step $\ke_\omega\rightsquigarrow\ke_{\bar{\omega}}$ will be used
in order to prove that $H^2 S^2\ke_\omega=0$ follows from $H^2
S^2\ke_{\bar{\omega}}=0$ even if $h^1\ke_{\bar{\omega}}(1)=1$.

As before we use the following notation. $\ke$ will denote a $5$--instanton and 
$\bar{\ke}$
the $(4,4)$--instanton obtained from a general $\xi\in H_5^\ast\cong
H^1\ke(-1)$ by the process $\omega\mapsto \bar{\omega}=\res_\xi \omega,$ such
that $\ke=\ke_\omega$ and $\bar{\ke}=\ke_{\bar{\omega}}$, see \ref{exxi}.

By Proposition \ref{exxi} we may assume that there is a second element
$\eta\in H^1\ke(-1)$, linear independent of $\xi$, such that the multiplication
map $\langle \xi,\eta\rangle \otimes V^\ast\to H^1\ke$ is injective (in fact an
isomorphism). Let $\bar{H}=\Ker(\xi)\subset H_5$ and let $\bar{\eta}\in
\bar{H}_4^\ast$ be the image of $\eta$. If 
$\bar{\omega}=\res_\xi\omega$, then $\bar{\omega}\in M(4,4)_{\bar{\eta}}$, and
we have $H^2(S^2\bar{\ke})=0$, see \ref{vanind} in case $n=4$. Now $\ke$
appears as the cohomology of the monad
\[
0\to \ko_\P (-1)\to\bar{\ke}\to\ko_\P (1)\to 0
\]
resulting from $\xi$. Then $S^2\ke$ appears as the cohomology of the derived
monad
\[
0\to \bar{\ke}(-1)\to S^2\bar{\ke}\oplus\ko\to\bar{\ke}(1)\to 0,
\]
see proof of \ref{crism} and \ref{vanind}. Because of the vanishing conditions
for instantons, we obtain an exact sequence
\[
H^1(S^2\bar{\ke})\to H^1\bar{\ke}(1)\to H^2(S^2\ke)\to H^2(S^2\bar{\ke})=0.
\]
By Proposition \ref{boundh1} we may assume that $h^1\bar{\ke}(1)\leq 1$. In
what follows, we shall, practically, describe the kernel of the map
$H^2\bigl(\bar{\ke}(1)\bigr)^\ast \to H^1 (S^2\bar{\ke})^\ast$ and
show that it is zero. Then it follows that $H^2S^2\ke=0$. In order to do that,
we first describe the dual of $H^1\bar{\ke}(1)$ in terms of its defining
tensor $\bar{\omega}$, which can be done for arbitrary $n$.  \vskip5mm

\begin{sub}\label{descrh1} {\bf Lemma}: Let $H$ be an $n$--dimensional vector
  space, let $\omega\in S^2H^\ast\otimes \Lambda^2 V^\ast$ be non--degenerate
  and let $\ke=\ke_\omega$, and interpret the elements of $H\otimes S^2 V$ as
  linear maps $H^\ast\otimes V^\ast\to V$ which are selfdual with respect to
  $V$. Then 
  \begin{enumerate}
  \item [(a)] $H^1(\ke(1))^\ast\cong\{\gamma\in H\otimes S^2 V\ |\ \gamma\circ
    \omega=0\}$
\item [(b)] If $P=\P W$ is a plane in $\P_3,\, W\subset V$, then
  $H^1(\ke_P(1))^\ast\cong \{\gamma\in H\otimes S^2 V\ |\ \gamma\circ
  \omega=0\ \text{ and } \Im\ \gamma\subset W\}$
  \end{enumerate}
\end{sub}
\vskip5mm

\begin{proof}
a)\; From the selfdual monad associated to $\omega$, one derives the exact
sequence
\[
N \otimes V^\ast \to H^\ast \otimes S^2 V^\ast \to H^1 \ke(1) \to 0
\]
and its dual
\[
0 \to H^1\bigl(\ke(1)\bigr)^\ast \to H \otimes S^2 V \to N^\ast \otimes V.
\]

Interpreting the elements of $N^\ast \otimes V$ as linear maps $N \to V$ and
the elements of $H \otimes S^2 V$ as linear maps $H^\ast \otimes V^\ast
\xrightarrow{\gamma} V \cong H_1 \otimes V$, the homomorphism in the last
sequence can be described by $\gamma \mapsto \gamma \circ u$, where $u : N
\hookrightarrow H^\ast \otimes V^\ast$.

Since $N = \Im(\omega)$, it follows that $\gamma \circ u = 0$ if and only if
$\gamma \circ \omega = 0$.

b)\; If $\gamma \in H \otimes S^2 V$ and $\Im(\gamma) \subset W$, then $\gamma
\in H \otimes S^2 W$, i.e.\ $\gamma$ factorizes as
\[
\xymatrix{
H^\ast \otimes V^\ast \ar[r]^>>>>>\gamma \ar@{>>}[d] & V\\
H^\ast \otimes W^\ast \ar[r]^>>>>>{\gamma_0} & W\ar@{>->}[u]
}
\]
for some $\gamma_0$ which is selfdual with respect to $W$.   Now restrict the
monad of $\omega$ to $P$ and use the same argument as in a).
\end{proof}
\vskip5mm

\begin{sub}\label{prsm} {\bf Proposition:}
Any $(5,2)$--instanton $\ke$  on $\P_3$ satisfies $H^2 S^2 \ke = 0$.
\end{sub}
\vskip5mm

\begin{proof}
a)\; Let $\omega \in M(5,2)$ define $\ke = \ke_\omega$.   As before, we
denote by $\bar{\ke}$ the bundle $\ke_{\bar{\omega}}$ defined by $\bar{\omega}
= \res_\xi(\omega)$ for a general $\xi \in H^\ast$.   We may assume that $h^1
\bar{\ke}(1) \le 1$ by Proposition \ref{boundh1}.   As mentioned above, we may
also assume that $H^2 S^2 \bar{\ke} = 0$.

b)\; Let $H_1 \oplus \bar{H} = H$ be a decomposition defined by $\xi$ with
$\bar{H} = \Ker(\xi)$.   Then
\[
(H_1 \oplus \bar{H}) \otimes V \xrightarrow{\omega} (H^\ast_1 \oplus
\bar{H}^\ast) \otimes V^\ast
\]
decomposes as
\[
\omega = \left(
  \begin{array}{c|c}
e \circ \bar{\varphi} \circ e^\ast & e \circ \bar{\varphi} \circ \bar{u}^\ast
\\ \hline
\bar{u} \circ \bar{\varphi} \circ e^\ast & \bar{\omega}
  \end{array}\right)
\]
where $\bar{N} \xrightarrow{e} H_1^\ast \otimes V^\ast$ is a linear map and
where
\[
\xymatrix{
\bar{H} \otimes V \ar[r]^>>>>{\bar{\omega}}\ar[d]^>>>>{\bar{u}^\ast} &
\bar{H}^\ast \otimes V^\ast\\ 
\bar{N} \ar[r]^>>>>>>>>{\bar{\varphi}}_>>>>>>>>\cong & \bar{N} \ar[u]_>>>>{\bar{u}}
}
\]
is the decomposition of $\bar{\omega} = \bar{u} \circ \bar{\varphi} \circ
\bar{u}^\ast$.   Note that $\omega$ is symmetric, resp.\ skew with respect to
$H$, resp.\ $V$.   Similarly, any $\sigma \in \Lambda^2 H \otimes S^2 V$ can
be written as a linear operator
\[
(H_1^\ast \oplus \bar{H}^\ast) \otimes V^\ast \xrightarrow{\sigma} (H_1 \oplus
\bar{H}) \otimes V
\]
and decomposed as
\[
\sigma = \left(
  \begin{array}{c|c}
0 & \gamma\\ \hline
-\gamma^\ast & \bar{\sigma}
  \end{array}\right)
\]
with $\gamma \in \bar{H} \otimes S^2 V$ because $\sigma$ is skew with respect
to $H$. 

c)\; We are now using the exact sequence
\[
0 \to H^2(S^2 \ke)^\ast \to \Lambda^2 H \otimes S^2 V \to N^\ast \otimes H
\otimes V
\]
of Remark \ref{remsm}, the homomorphism being $\sigma \mapsto \sigma \circ
\omega$.   By the above decomposition the condition $\sigma \circ \omega = 0$
is equivalent to the two conditions
\[
\gamma \circ \bar{u} = 0 \text{ and } \gamma^\ast \circ e = \bar{\sigma} \circ
\bar{u}.
\]
In order to show that $\sigma \circ \omega = 0$ implies $\sigma = 0$, it is
now sufficient to prove that the two conditions on $\gamma$ imply that $\gamma
= 0$.   For then, also, $\bar{\sigma} \circ \bar{u} = 0$ or $\bar{\sigma}
\circ \bar{\omega} = 0$, and then $\bar{\sigma} = 0$, because $H^2S^2
\bar{\ke} = 0$, using the same sequence for $\bar{\ke}$.

d)\; The first condition $\gamma \circ \bar{u} = 0$ or $\gamma \circ
\bar{\omega} = 0$ means that
\[
\gamma \in H^1\bigl(\bar{\ke}(1)\bigr)^\ast \subset \bar{H} \otimes S^2 V
\]
by Lemma \ref{descrh1}, (a).   If $H^1 \bar{\ke}(1) = 0$, there is nothing to
prove.   So we may assume that $\gamma \not= 0$.   By Lemma \ref{vanpl} we may
assume that $H^1 \bar{\ke}_P(1) = 0$ for any plane $P = \P W$ in $\P V$.
Then, by Lemma \ref{descrh1}, (b), $\Im(\gamma) \not\subset W$ for any
3--dimensional subspace $W \subset V$.   Therefore, 
\[
\bar{H}^\ast \otimes V^\ast \xrightarrow{\gamma} H_1 \otimes V
\]
is surjective.   Because $\gamma \circ \bar{u} = 0$ and $\dim \bar{N} = 12$,
we obtain the exact diagram
\[
\xymatrix{
0 \ar[r] & \bar{N} \ar[r] & \bar{H}^\ast \otimes V^\ast \ar[r]^>>>>\gamma &
H_1 \otimes V \ar[r] & 0\\
0 \ar[r] & N \ar[r]\ar[u]_>>>>\cong & H^\ast \otimes V^\ast \ar[r] \ar@{>>}[u]
& Q \ar[r]\ar@{>>}[u] & 0\\
& & \xi \otimes V^\ast \ar[u] & & \\
& & 0. \ar[u] & &
}
\]
Note that the induced projection induces an isomorphism between $N$ and
$\bar{N}$, because $\omega$ and $\bar{\omega}$ have the same rank 12.   It
follows from the decomposition of $\omega$ that any vector of $N$ can be
written uniquely as
\[
e (\theta) + \theta
\]
with $\theta \in \bar{N} \subset \bar{H}^\ast \otimes V^\ast$, with the
isomorphism given by $e(\theta) + \theta \leftrightarrow \theta$, and $\Im(e)
\subset \xi \otimes V^\ast$. 

e)\; We are now going to show that the condition $\gamma^\ast \circ e =
\bar{\sigma} \circ \bar{u}$ implies a contradiction.   For that we choose a
2--dimensional subspace $K \subset H^\ast$ containing $\xi$, using Lemma
\ref{indec}, such that $N \cap (K \otimes V^\ast)$ contains no non--zero
decomposable vector.   Let $\xi, \xi^\prime$ be the basis of $K$, $\xi^\prime
\in \bar{H}^\ast$.   Then there are two independent linear forms $z, z^\prime
\in V^\ast$ such that $N \cap (K \otimes V^\ast)$ is spanned by
\[
\xi \otimes z + \xi^\prime \otimes z^\prime.
\]
Because $\xi \otimes z \mapsto 0$ under $N \to \bar{N}$, we have with $\theta =
\xi^\prime \otimes z^\prime$ that
\[
\xi \otimes z + \xi^\prime \otimes z^\prime = e (\xi^\prime \otimes z^\prime)
+ \xi^\prime \otimes z^\prime
\]
or $e (\xi^\prime \otimes z^\prime) = \xi \otimes z$.   Therefore,
$\gamma^\ast (\xi \otimes z) = \bar{\sigma} (\xi^\prime \otimes z^\prime)$.

Because $\bar{\sigma} \in \Lambda^2 \bar{H} \otimes S^2 V$ is skew with
respect to $\bar{H}$, it follows $\bar{\sigma}(\eta \otimes w) \in \Ker(\eta)
\otimes V$ under $\bar{H}^\ast \otimes V^\ast \xrightarrow{\bar{\sigma}}
\bar{H} \otimes V$ for any $\eta \in \bar{H}^\ast$ and any $w \in V^\ast$.   It
follows that $(\xi^\prime \otimes \id) \circ \gamma^\ast (\xi^\prime \otimes
z) = 0$, and dually, that $\gamma (\xi^\prime \otimes z) = 0$, by considering
the diagrams
\[
\begin{array}{cc}
\xymatrix{
H^\ast_1 \otimes V^\ast \ar[r]^>>>>{\gamma^\ast} & \bar{H} \otimes V
\ar[d]^>>>>{\xi^\prime \otimes \id} \\
\langle \xi \otimes z\rangle \ar[u] \ar[r]^>>>>>>0 & V
} &\qquad\qquad
\xymatrix{
H_1 \otimes V\ar[d]^>>>>{\xi \otimes z} & \bar{H}^\ast \otimes 
V^\ast\ar[l]_>>>>{\gamma^\ast}\\
k & \xi^\prime \otimes V^\ast \ar[l]_>>>>>>>0 \ar[u]
}
\end{array}
\]
and using that $\gamma$ is symmetric with respect to $V$.

Now $\xi^\prime \otimes z \in \bar{N}$ and, because $e$ has its image in $\xi
\otimes V^\ast$, there is a form $z^{\prime\prime} \in V^\ast$ with
$e(\xi^\prime \otimes z) = \xi \otimes z^{\prime\prime}$.   Then
\[
\xi \otimes z^{\prime\prime} + \xi^\prime \otimes z \in N \cap (K \otimes
V^\ast)
\]
is a second vector which is independent of $\xi \otimes z + \xi^\prime \otimes
z^\prime$, a contradiction to the choice of $K$.   Therefore, the condition
$\gamma^\ast \circ e = \bar{\sigma} \circ \bar{u}$ for $\gamma \in
  H^1\bigl(\bar{\ke}(1)\bigr)^\ast$ implies $\gamma = 0$, which proves
  Proposition \ref{prsm}. 
\end{proof}

\end{document}